\documentclass{amsart}

\pdfoutput=1

\usepackage[margin=1.25in]{geometry}

\usepackage[utf8]{inputenc}
\usepackage[T1]{fontenc}

\usepackage[backend=biber,style=alphabetic]{biblatex}
\usepackage{amsmath}
\usepackage{amsfonts}
\usepackage{amssymb}
\usepackage{amsthm}
\usepackage{tikz-cd}
\usepackage{hyperref}
\usepackage{cleveref}
\usepackage{enumerate}
\usepackage{xcolor}
\usepackage{mathtools}
\usepackage{capt-of}

\usepackage{spectralsequences}
\definecolor{pastelblue}{rgb}{0.25,0.54,0.8}
\definecolor{pastelred}{rgb}{0.8,0.2,0.3}
\definecolor{pastelgreen}{rgb}{0.2,0.8,0.4}
\definecolor{pastelorange}{rgb}{0.9,0.4,0.1}

\usepackage{multicol}

\newtheorem{thm}{Theorem}
\newtheorem{prop}[thm]{Proposition}
\newtheorem{cor}[thm]{Corollary}
\newtheorem{conj}[thm]{Conjecture}

\newtheorem*{thm*}{Theorem}
\newtheorem*{prop*}{Proposition}

\theoremstyle{definition}
\newtheorem{defn}[thm]{Definition}
\newtheorem{rem}[thm]{Remark}

\newcommand{\ext}{\text{Ext}}

\newcommand{\sq}{\text{Sq}^0}

\title{The $\mathbb{C}$-Motivic Adams Spectral Sequence for $s\leq5$}
\author{Jordan Benson}
\date{\textit{Last Revised} \today}

\begin{document}

\begin{abstract}
	We determine the $\tau^n$-torsion in the first 5 lines of the $E_2$ page of the $\mathbb{C}$-motivic Adams spectral sequence using the techniques of Burklund-Xu. In particular, every element in this range is either $\tau^1$-torsion or $\tau$-free. We also show that $\tau^n$-torsion elements can appear only in Adams filtration at least $2n+2$ and give further evidence of a possible $3n$ bound.
\end{abstract}

\maketitle

\tableofcontents

\section{Introduction}

The foremost computational problem in stable homotopy theory is to determine the stable homotopy groups of spheres, which fit into a graded ring. We will refer to the graded components of this ring as stable stems.

%In principle, if the entire multiplicative structure of this ring is known, including all matric Toda brackets of elements, then one can compute the stable homotopy groups of any cellular object using the Atiyah-Hirzebruch spectral sequence. Knowledge of stable stems can also be applied to low-dimensional topology. For example, one can classify boundaries of parallelizable manifolds, as in \textcite{HHR}, and determine the existence of exotic smooth structures, as in \textcite{WangXu61}.

Beyond applications to problems in homotopy theory itself, knowledge of stable stems can also be used to solve problems in algebra and geometry. For example, there is a well-known correspondence between finite-dimensional division algebras over $\mathbb{R}$ and the Hopf invariant one elements in the 2-primary stable stems. \textcite{Adams} computed differentials on the Hopf invariant one elements and used this correspondence to prove that $\mathbb{R}$, $\mathbb{C}$, $\mathbb{H}$, and $\mathbb{O}$ are the only such algebras.

In seemingly unrelated work, \textcite{Browder} reduced the existence of manifolds with Kervaire invariant one to the existence of elements in the 2-primary stable stems that behave algebraicly as squares of the Hopf invariant one elements . Using this relationship and techniques from equivariant and chromatic stable homotopy theory, \textcite{HHR} showed that no such manifold exists beyond dimension 126. This is closely related to the work of \textcite{Kervaire-Milnor} on the classification of exotic smooth structures on spheres. In the 61-dimensional case, \textcite{Wang-Xu} used computations in the Adams spectral sequence to show that the 61-sphere has a unique smooth structure, the last possible case in odd dimension.

In order to study the 2-primary component of the stable homotopy groups of spheres, Adams introduced a spectral sequence of the form
\[
	E_2^{s,t} = \text{Ext}_\mathcal{A}^{s,t}(\mathbb{F}_2,\mathbb{F}_2)\Rightarrow\pi_{t-s}(S^0)\otimes\mathbb{Z}_{2},
	\]
where $S^0$ denotes the sphere spectrum, $\text{Ext}_\mathcal{A}^{\ast,\ast}$ denotes Ext in comodules over the dual Steenrod algebra, and $\mathbb{Z}_2$ denotes the 2-adic integers. Henceforth, we work implicitly at the prime 2.

Recently, \textcite{GWX} and \textcite{IWX} have had major success in computing stable stems by exploiting a curious property of the cellular $\mathbb{C}$-motivic category. To summarize that approach, recall that there is a map $\tau:\widehat{S^{0,-1}}\to \widehat{S^{0,0}}$, where $\widehat{S^{\ast,\ast}}$ denotes a 2-completed $\mathbb{C}$-motivic sphere. Importantly, the map $\tau$ exists because we are working at a prime, rather than integrally. This map $\tau$ has the property that its cofiber $C\tau$ is an $\mathbb{E}_\infty$-algebra, by \textcite{Gheorghe} (see also \cite{GWX} and \textcite{Piotr}), and furthermore satisfies
\[
	\pi_{\ast,\ast}C\tau\cong\text{Ext}^{\ast,\ast}_{BP_\ast BP}(BP_\ast,BP_\ast),
	\]
by \textcite{isaksen14}. The right-hand side is seen to be the $E_2$ page of the 2-localized Adams-Novikov spectral sequence, whose abutment is the 2-localized stable stems. Furthermore, upon inverting the element $\tau$, there is an isomorphism
\[
	\pi_{\ast,\ast}\widehat{S^{0,0}}\otimes_{\mathbb{Z}_2[\tau]}\mathbb{Z}_2[\tau^{-1}]\cong\pi_\ast S^0\otimes_{\mathbb{Z}_2}\mathbb{Z}_2[\tau^{\pm1}]
\]
by \textcite{DI}. Then we see that computing motivic stable stems is sufficient for computing classical stable stems. As a result, one may wish to set up an analog of the Adams spectral sequence in this motivic setting and use it to study motivic stable stems. In fact, the form of this spectral sequence is very similar to the classical one:
\[
	\text{Ext}^{\ast,\ast,\ast}_{\mathcal{A}^\text{mot}}(\mathbb{F}_2[\tau],\mathbb{F}_2[\tau])\Rightarrow\pi_{\ast,\ast}(\widehat{S^{0,0}})
\]
Here we use $\mathcal{A}^\text{mot}$ to denote the $\mathbb{C}$-motivic dual Steenrod algebra. Now, since $C\tau$ is a 2-cell complex, we have a long exact sequence
\[
	\dots\xrightarrow{\tau}\pi_{\ast,\ast}\widehat{S^{0,0}}\to\pi_{\ast,\ast}C\tau\to\pi_{\ast,\ast}\widehat{S^{1,-1}}\xrightarrow{\tau}\dots
\]
establishing a close connection between the motivic stable stems and the Adams-Novikov $E_2$ page. Using this long exact sequence and the motivic Adams spectral sequence, one can set up a motivic analog of the classical Miller square \cite{Miller}. By comparing partial information in two different spectral sequences, one can obtain a number of differentials not apparent in either individual spectral sequence, allowing for a substantially larger range of stable stems to be computed.

In order to use this approach, one needs as input some information about the $\mathbb{C}$-motivic Adams spectral sequence. This is the content of our main results. Before stating them, note that Proposition 3.5 of \cite{DI} provides an identification of the $\tau$-free part of the $\mathbb{C}$-motivic Adams $E_2$ page with the classical Adams $E_2$ page, similar to the identification in homotopy described above:

\begin{prop*}[\cite{DI}]
	There is a ring isomorphism
	\[
		\ext^{\ast,\ast,\ast}_{\mathcal{A}^\text{mot}}(\mathbb{F}_2[\tau],\mathbb{F}_2[\tau])\otimes_{\mathbb{F}_2[\tau]}\mathbb{F}_2[\tau^{\pm1}]\xrightarrow{\cong}\ext^{\ast,\ast}_\mathcal{A}(\mathbb{F}_2,\mathbb{F}_2)\otimes_{\mathbb{F}_2}\mathbb{F}_2[\tau^{\pm1}]
	\]
\end{prop*}

Then the classical Adams $E_2$ page accounts for the $\tau$-free part of the $\mathbb{C}$-motivic one, so it remains to identify the $\tau^n$-torsion elements in the range of interest:

\begin{thm}\label{thm:5-line}
	As an $\mathbb{F}_2[\tau]$-module, we have
	\[
		\text{Ext}^{\leq5,\ast,\ast}_{\mathcal{A}^\text{mot}}(\mathbb{F}_2[\tau],\mathbb{F}_2[\tau])\cong F\oplus T,
	\]
	where
	\[
		F\cong\text{Ext}^{\leq5,\ast}_\mathcal{A}(\mathbb{F}_2,\mathbb{F}_2)[\tau]
	\]
	and $T$ is an $\mathbb{F}_2$-vector space generated by
	\[
		\begin{array}{c c c c}
			h_1^4\text{ in }(4,4,4) &h_1^2c_0\text{ in }(10,5,7) &h_3g\text{ in }(27,5,16) &h_1^4h_{j+5}\text{ in }(3+2^{j+5},5,4+2^{j+4})
		\end{array}
	\]
	in the stated $(t-s,s,w)$-tridegrees, where $w$ is the $\mathbb{C}$-motivic weight.
\end{thm}

%\begin{thm}\label{thm:5-line}
%	Through the 5-line, the $\mathbb{C}$-motivic Adams $E_2$ page has as $\mathbb{F}_2[\tau]$-generators all the $\mathbb{F}_2$-generators of the classical Adams $E_2$ page, along with the $\tau$-torsion classes (bidegree $(t-s,s)$ listed)
%
%	\[
%		\begin{tabular}{c c c c}
%			$h_1^4$ in $(4,4)$ &$h_1^2c_0$ in $(10,5)$ &$h_3g$ in $(27,5)$ &$h_1^4h_{j+5}$ in $(2^{j+5}+3,5)$
%		\end{tabular}
%	\]
%
%	\noindent for $j\geq0$. Furthermore, the $\mathbb{F}_2$-module relations on the classical Adams $E_2$ page determine all the $\mathbb{F}_2[\tau]$-module relations on the $\mathbb{C}$-motivic Adams $E_2$ page, adjusting for motivic weight. In particular, the $\tau$-torsion classes participate in no relations not required by the ring structure.
%\end{thm}
%
%\begin{rem}
%	This theorem does not address the motivic weights of all generators. By ``adjusting for motivic weight," we mean adding appropriate powers of $\tau$ to classical relations to ensure that motivic weight is respected. By a relation ``required by the ring structure," we mean, for example, that $h_1\cdot h_1^3 = h_1^4$ because $h_1$ is an algebra generator. That the $\tau$-torsion classes participate in no new relations beyond these follows purely from a degree argument.
%\end{rem}

\begin{rem}\label{rem:weight}
	We make no claim about the motivic weights of the elements in $F$, but if one does know them, then one can easily deduce relations among these elements. For example, if one knows that $h_0$ has weight 0 while $h_{j+1}$ has weight $2^j$ for $j\geq0$, then the classical relation
	\[
		h_1^3 = h_0^2h_2
	\]
	no longer holds on the $\mathbb{C}$-motivic Adams $E_2$ page. Instead, the weight of the left-hand side needs to be lowered by 1 in order to match the weight of the right-hand side:
	\[
		\tau h_1^3 = h_0^2h_2
	\]
\end{rem}

\begin{rem}
	This theorem implies, for example, that the $\tau$-power torsion appearing through the 5-line in the charts of \textcite{charts} is all that exists in the $\mathbb{C}$-motivic Adams $E_2$ page, except the various $h_1^4h_{j+1}$ which are out of the computed range of stems in those charts. In particular, the highest power of $\tau$-torsion appearing through the 5-line is $\tau^1$-torsion.
\end{rem}

As a corollary of \Cref{thm:5-line}, we can deduce the $\mathbb{C}$-motivic Hopf differentials, following the argument found in \textcite{Balderrama} in the $\mathbb{R}$-motivic case and \textcite{Wang} in the classical case:

\begin{cor}
	There are differentials
	\[
		d_2(h_j) = h_0h_{j-1}^2
	\]
	in the $\mathbb{C}$-motivic Adams spectral sequence for all $j\geq4$.
\end{cor}
\begin{proof}
	We first note that graded-commutativity of the classical stable stems $\pi_\ast S^0$ force
	\[
		2\sigma^2 = 0,
	\]
	where $\sigma\in\pi_7S^0$ is a stable Hopf map. Since 2 is detected on the Adams $E_2$ page by $h_0$ and $\sigma$ by $h_3$, the class $h_0h_3^2\neq0$ on the Adams $E_2$ page cannot survive to homotopy. If it did, then it would detect $2\sigma^2$.

	By \Cref{thm:5-line}, there is a $\tau$-free generator on the $\mathbb{C}$-motivic Adams $E_2$ page which is a lift of $h_0h_3^2$. Abusing notation by denoting this element by $h_0h_3^2$, there must be a motivic Adams differential received by $h_0h_3^2$. Indeed, this element is a product of permanent cycles, so the Leibniz rule forces $h_0h_3^2$ itself to be a permanent cycle. The only element which can possibly support a differential hitting $h_0h_3^2$ is $h_4$, giving the differential
	\[
		d_2(h_4) = h_0h_3^2
	\]
	on the $\mathbb{C}$-motivic Adams $E_2$ page.

	Now suppose that we have obtained a differential
	\[
		d_2(h_j) = h_0h_{j-1}^2
	\]
	for some $j\geq4$. By the $\mathbb{C}$-motivic relation
	\[
		h_jh_{j+1} = 0
	\]
	following from the analogous classical relation, we have
	\[
		0 = d_2(h_jh_{j+1}) = h_0h_{j-1}^2h_{j+1}+h_jd_2(h_{j+1})
	\]
	by the Leibniz rule. Since $j\geq4$ by assumption, the comments in \Cref{rem:weight} indicate that the classical relation
	\[
		h_{j-1}^2h_{j+1} = h_j^3
	\]
	holds on the $\mathbb{C}$-motivic Adams $E_2$ page without modification by powers of $\tau$. As a result, we find that
	\[
		h_0h_j^3 = h_jd_2(h_{j+1})
	\]
	with $d_2(h_{j+1})$ lying in a tridegree with sole generator $h_0h_j^2$. Then it suffices to show that $h_0h_j^3\neq0$. Since this holds classically, \Cref{thm:5-line} implies that it holds on the $\mathbb{C}$-motivic Adams $E_2$ page. Therefore, we have the differential
	\[
		d_2(h_{j+1}) = h_0h_j^2,
	\]
	as desired.
\end{proof}
%\begin{proof}
%	The element
%	\[
%		h_3\in\ext^{1,8,4}_{\mathcal{A}^\text{mot}}(\mathbb{F}_2[\tau],\mathbb{F}_2[\tau])
%	\]
%	survives to the bigraded homotopy of the sphere so detects homotopy classes. One of these classes Betti realizes to the classical element
%	\[
%		\sigma\in\pi_7S^0,
%	\]
%	whose square $\sigma^2$ must be 2-torsion since $\sigma$ lives in an odd stem, by graded-commutativity of $\pi_\ast S^0$. If $h_0h_3^2$ were to survive the $\mathbb{C}$-motivic Adams spectral sequence, then it would detect a $\mathbb{C}$-motivic homotopy class Betti realizing to $2\sigma^2=0$. As a result, this class must participate in a $\mathbb{C}$-motivic Adams differential.
%
%	Since $h_0h_3^2$ is a product of permanent cycles, it must also be a permanent cycle, by the Leibniz rule. Then this class must receive a differential. The only possibility is
%	\[
%		d_2(h_4) = h_0h_3^2
%	\]
%	for degree reasons. Now suppose that
%	\[
%		d_2(h_j) = h_0h_{j-1}^2
%	\]
%	for some $j\geq4$. By the relation
%	\[
%		h_jh_{j+1} = 0
%	\]
%	on the $\mathbb{C}$-motivic Adams $E_2$ page, the Leibniz rule gives
%	\[
%		0 = d_2(h_jh_{j+1}) = d_2(h_j)h_{j+1}+h_jd_2(h_{j+1}),
%	\]
%	which gives
%	\[
%		h_0h_{j-1}^2h_{j+1}+h_jd_2(h_{j+1})
%	\]
%	under the inductive hypothesis.
%\end{proof}

Turning now to the more general pattern of $\tau^n$-torsion on the $\mathbb{C}$-motivic Adams $E_2$ page, in the computed range of stems in \cite{charts}, it seems that higher $\tau^n$-torsion appears only in higher filtration. In fact this holds for all stems:

\begin{thm}\label{thm:generalbound}
	If $x$ is a class on the $\mathbb{C}$-motivic Adams $E_2$ page and $\tau^nx=0$ while $\tau^{n-1}x\neq0$, then $x$ has Adams filtration at least $2n+2$.
\end{thm}

This bound is not tight, and based on the charts of \cite{charts}, it seems reasonable to expect that the $2n+2$ above could be replaced by a better bound:

\begin{conj}
	If $x$ is a class on the $\mathbb{C}$-motivic Adams $E_2$ page and $\tau^nx=0$ while $\tau^{n-1}x\neq0$, then $x$ has Adams filtration at least $3n$.
\end{conj}

If it were to hold in general, this bound would be tight, by the existence of the $\tau^3$-torsion element $h_2g^2$ in filtration 9. By \Cref{thm:5-line}, the $3n$ bound holds through the 5-line. In fact, we can prove a better bound for $\tau^2$-torsion elements, even though we do not have enough information to compute the entire 6-line by our methods:

\begin{thm}\label{thm:6-line}
	Through the 6-line, all $\mathbb{C}$-motivic Adams $E_2$ classes are either $\tau$-free or $\tau^1$-torsion.
\end{thm}

In particular, all $\mathbb{C}$-motivic Adams $E_2$ classes which are $\tau^2$-torsion but not $\tau$-torsion can appear only in the 7-line and up, giving further empirical evidence for the conjectured $3n$ bound.

By \Cref{thm:generalbound}, there is no $\tau^n$-torsion appearing until at least the 4-line. The relations in this range are, up to powers of $\tau$ to account for motivic weight, exactly those found on the classical Adams $E_2$ page. Informally, one may think of the motivic Adams $E_2$ page as being the ``$\tau$-adjoined'' version of the classical one in this range. However, by low-dimensional calculations, there is a $\tau$-torsion element, namely $h_1^4$, living on the 4-line. This is closely related to the discussion in \Cref{rem:weight}.

Classically, \textcite{Lin} and his student \textcite{Chen} computed the 4- and 5-lines of the Adams $E_2$ page using the Lambda algebra. In this method, complete information about the $s$-line is taken as input. Calculations are made in low internal degrees on the $(s+1)$-line using the Lambda algebra and this input, and then one uses this to compute one degree higher by relating the Adams $E_2$ page for $\mathbb{R}P^\infty_{-\infty}$ to that for the sphere spectrum, which is critical to the proof of Lin's theorem \cite{Lin-KP}. The combination of these steps allows one to use low-degree calculations inductively to move along the $(s+1)$-line and work out all elements there. The result of Lin and Chen's work is summarized in \cite{Chen} and provides an $\mathbb{F}_2$-module presentation for the Adams $E_2$ page for Adams filtration $s\leq5$.

The approach we follow is the one laid out in \textcite{BX} and uses complete information about the classical Adams $E_2$ page through the $s$-line to work out the $\mathbb{C}$-motivic Adams $E_2$ page through the $s$-line. This approach can be summarized in the following way, where algAHSS denotes the so-called \textit{algebraic Atiyah-Hirzebruch spectral sequence} of \cite{BX} and CESS denotes the usual Cartan-Eilenberg spectral sequence:
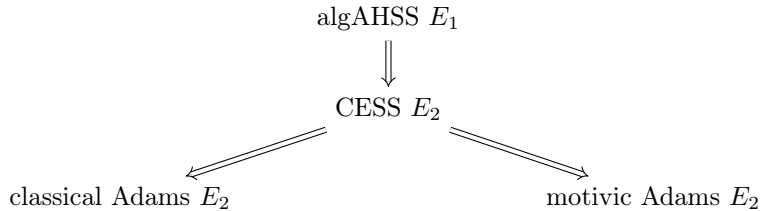
\[
	\begin{tikzcd}
		&\text{algAHSS }E_1\arrow[Rightarrow]{d} & \\
		&\text{CESS }E_2\arrow[Rightarrow]{dl} \arrow[Rightarrow]{dr} & \\
		\text{classical Adams }E_2 & &\text{motivic Adams }E_2
	\end{tikzcd}
\]
This diagram will be explained in the following sections. In \Cref{sec:CESS} we describe the main properties of the CESS which we will need for our computations. In particular, we show that there are only a few differential patterns we need to consider, based on sparsity of the CESS. We also prove \Cref{thm:generalbound} by appealing to a rigidity theorem of \cite{BX}. In \Cref{sec:algAHSS} we describe in detail the algAHSS, the spectral sequence in which we will be carrying out the bulk of our computations. We show how one computes differentials in this spectral sequence and prove further simplifying properties beyond those proven in \cite{BX}.

Once we have introduced the spectral sequences above, we prove the main results in \Cref{sec:main}. The proofs proceed by first establishing which algAHSS elements survive to the CESS, of those which can possibly contribute $\tau^n$-torsion elements in the desired range. We then show that these elements either survive to the classical Adams $E_2$ page or participate in CESS differentials by appealing to the low-dimensional calculations in \cite{IWX} and various degree arguments.

\section*{Acknowledgments}

The author would like to thank his advisor, Zhouli Xu, for suggesting this project and giving feedback on early drafts of this paper. Among other things, he has passed along the sage advice of Mark Mahowald that one should ``get to know homotopy elements by their first names,'' and this project certainly went a long way toward that end. The author would also like to thank Maxwell Johnson, Shangjie Zhang, and Dan Isaksen for feedback on an early draft of this paper. The author would especially like to thank Eva Belmont for detailed feedback and innumerable meetings while this project was underway.

\section{The Cartan-Eilenberg Spectral Sequence}\label{sec:CESS}

We begin this section by summarizing some of the work in \cite{BX} and describing the differential pattern in the Cartan-Eilenberg spectral sequence. Recall that there is a Hopf algebra extension
\[
	\mathcal{P}\hookrightarrow\mathcal{A}\to\mathcal{Q},
\]
where $\mathcal{A}$ is the mod-2 dual Steenrod algebra and $\mathcal{P}$ is the ``doubled-up'' Hopf subalgebra
\[
	\mathcal{P} = \mathbb{F}_2[\xi_1^2,\xi_2^2,\dots]\subset\mathcal{A}
\]
with the induced structure maps. Then
\[
	\mathcal{Q} = \mathbb{F}_2[\xi_1,\xi_2,\dots]/(\xi_1^2,\xi_2^2,\dots)
\]
is the exterior algebra on the generators of $\mathcal{A}$. The form of the resulting Cartan-Eilenberg spectral sequence, which we refer to as the CESS, is
\[
	E_2^{s,k,t} = \text{Ext}_\mathcal{P}^s(\mathbb{F}_2,\text{Ext}_{\mathcal{Q}}^{k,t}(\mathbb{F}_2,\mathbb{F}_2))\Rightarrow\text{Ext}_\mathcal{A}^{s+k,t}(\mathbb{F}_2,\mathbb{F}_2),
\]
where the abutment is seen to be the classical Adams $E_2$ page.

In order to compute the $\mathbb{C}$-motivic Adams $E_2$ page, one can consider the analogous Hopf algebra extension
\[
	\mathcal{P}^\text{mot}\hookrightarrow\mathcal{A}^\text{mot}\to\mathcal{Q}^\text{mot},
\]
where $\mathcal{A}^\text{mot}$ is the $\mathbb{C}$-motivic mod-2 dual Steenrod algebra and $\mathcal{P}^\text{mot}$ and $\mathcal{Q}^\text{mot}$ are defined similarly to their classical counterparts.

The form of the resulting motivic Cartan-Eilenberg spectral sequence, which we will refer to also as the CESS, is
\[
	E_2^{s,k,t} = \text{Ext}_{\mathcal{P}^\text{mot}}^s(\mathbb{F}_2[\tau],\text{Ext}_{\mathbb{Q}^\text{mot}}^{k,t}(\mathbb{F}_2[\tau],\mathbb{F}_2[\tau]))\Rightarrow\text{Ext}_{\mathcal{A}^\text{mot}}^{s+k,t}(\mathbb{F}_2[\tau],\mathbb{F}_2[\tau])
\]

By construction, the $\mathbb{F}_2$-module basis for the classical CESS $E_2$ page provides an $\mathbb{F}_2[\tau]$-module basis for the motivic one, so we need only compute the classical $E_2$ page and then freely adjoin $\tau$. This justifies our convention of referring to both sequences as CESS.

In both sequences, the differentials are of the form
\[
	d_r:E_r^{s,k,t}\to E_r^{s+r,k-r+1,t}
\]
and so behave like Adams $d_1$ differentials if one remembers that $s+k$ is the total homological degree or, equivalently, total Adams filtration. Because every element of $\mathcal{P}$ or $\mathcal{P}^\text{mot}$ has even internal degree, we obtain a sparsity result for the CESS:

\begin{prop}[\cite{BX}]
If $r$ is even, then all CESS $d_r$ differentials are trivial.
\end{prop}

\noindent Then we can write every CESS differential as $d_{2n+1}$. Now we recall the rigidity theorem of \cite{BX} relating the classical and $\mathbb{C}$-motivic versions of the CESS.

\begin{prop}[\cite{BX}]\label{thm:rigidity}
	After a re-grading, the motivic Cartan-Eilenberg spectral sequence for
	\[
		\text{Ext}^{\ast,\ast,\ast}_{\mathcal{A}^\text{mot}}(\mathbb{F}_2[\tau],\mathbb{F}_2[\tau])
	\]
	is isomorphic to a $\tau$-B\"ockstein spectral sequence.
\end{prop}

\noindent In particular, if
\[
	d_{2n+1}(a) = b
\]
in the classical CESS, then
\[
	d_{2n+1}(a) = \tau^n b
\]
in the $\mathbb{C}$-motivic CESS. In other words, the differentials of the classical and $\mathbb{C}$-motivic CESS ``determine'' each other. Furthermore, the length of a classical CESS differential determines the power of $\tau$-torsion appearing in the target bidegree on the $\mathbb{C}$-motivic Adams $E_2$ page.

By appealing to this rigidity theorem, we will prove our main results by computing certain tridegrees in the classical CESS, comparing them against the classical Adams $E_2$ page, and then working backward to deduce the relevant CESS differentials. We then translate these differentials into $\tau^n$-torsion information about the $\mathbb{C}$-motivic Adams $E_2$ page.

Now recall that, since $\mathcal{Q}$ is an exterior $\mathbb{F}_2$-algebra on $\xi_i$ for $i\geq1$, the ring
\[
	\ext^{\ast,\ast}_\mathcal{Q}(\mathbb{F}_2,\mathbb{F}_2)
\]
has the structure of a ring isomorphic to
\[
	\mathbb{F}_2[q_0,q_1,\dots],
\]
with $q_i$ corresponding to $\xi_{i+1}$. By definition, this ring naturally forms a comodule over the ring $\mathcal{P} = \mathbb{F}_2[\xi_1^2,\xi_2^2,\dots]$ described above. One can check (\cite{BX}) that this comodule structure is as follows:
\[
	\psi(q_n) = \sum\limits_{i=0}^n\xi_{n-i}^{2^{i+1}}\otimes q_i
\]

We will return to this comodule structure momentarily, but we are now ready to prove \Cref{thm:generalbound}, which we reproduce here for convenience:

\begin{thm*}
	If $x$ is a class on the $\mathbb{C}$-motivic Adams $E_2$ page and $\tau^nx=0$ while $\tau^{n-1}x\neq0$, then $x$ has Adams filtration at least $2n+2$.
\end{thm*}
\begin{proof}[Proof of \Cref{thm:generalbound}]
	If $x$ is as described, then $x$ is detected by some CESS element $b$ receiving a CESS differential
	\[
		d_{2n+1}(a) = \tau^n b
	\]
	from another CESS element $a$. Because a CESS $d_{2n+1}$ differential lowers CE filtration by $r-1$ and the lowest possible CE filtration is 0, the CE filtration of $a$ must be at least $2n$. Since the Adams filtration of $a$ is the sum of its CE filtration $k$ and its homological degree $s$, the Adams filtration of $a$ must be at least $2n$.

	In fact, the Adams filtration of $a$ must be at least $2n+1$ since the only CESS element in bidegree $(s,k) = (0,k)$ is $q_0^k$, which must survive the CESS to detect $h_0^k$ on the Adams $E_2$ page. Indeed, sparsity forces $q_0^k$ to detect $h_0^k$.

	To see that no other $q_I$ survive to the CESS $E_2$ page, we compute the $s=0$ graded component of the CESS $E_2$ page directly. Recall that this $E_2$ page is of the form
	\[
		\text{Ext}^s_\mathcal{P}(\mathbb{F}_2,\ext^{k,t}_\mathcal{Q}(\mathbb{F}_2,\mathbb{F}_2)),
	\]
	so computing the $s=0$ graded component is equivalent to computing maps
	\[
		\mathbb{F}_2\to\ext^{k,t}_\mathcal{Q}(\mathbb{F}_2,\mathbb{F}_2)
	\]
	of $\mathcal{P}$-comodules. Let $f$ be such a map. Since the source of $f$ is $\mathbb{F}_2$, its image is determined entirely by $f(1)$, which must satisfy
	\[
		\psi(f(1)) = (1\otimes f)(\psi(1)),
	\]
	where we abuse notation by writing $\psi$ for the $\mathcal{P}$-comodule structures on both $\mathbb{F}_2$ and $\ext^{k,t}_\mathcal{Q}(\mathbb{F}_2,\mathbb{F}_2)$. The right-hand side reduces to
	\[
		(1\otimes f)(\psi(1)) = (1\otimes f)(1\otimes 1) = 1\otimes f(1),
	\]
	so $f(1)$ must be a primitive in
	\[
		\ext^{k,t}_\mathcal{Q}(\mathbb{F}_2,\mathbb{F}_2)\cong\mathbb{F}_2[q_0,q_1,\dots],
	\]
	where we have made an identification with $\mathbb{F}_2[q_0,q_1,\dots]$ for clarity. By the comodule structure described above, the only primitives in this $\mathcal{P}$-comodule are of the form $q_0^k$ for $k\geq0$.

	Then no CESS element in bidegree $(s,2n) = (0,2n)$ can support a CESS differential, so $a$ must lie in bidegree $(s,2n) = (\geq1,2n)$, hence must lie in Adams filtration at least $2n+1$. Since $a$ is the source of the differential in question and CESS differentials raise Adams filtration by 1, the target $\tau^nb$ must have Adams filtration at least $2n+2$.
\end{proof}

In order to simplify all the calculations that follow, we use the following result of \cite{BX}, which says that $\ext_\mathcal{Q}^{\ast,\ast}(\mathbb{F}_2,\mathbb{F}_2)$ admits a convenient splitting as a $\mathcal{P}$-comodule:

\begin{prop}[\cite{BX}]\label{prop:comodule}
	The $\mathcal{P}$-comodule structure
	\[
		\psi(q_n) = \sum\limits_{i=0}^n\xi_{n-i}^{2^{i+1}}\otimes q_i
	\]
	on $\mathbb{F}_2[q_0,q_1,\dots]$, where we set $\xi_0=1$ by convention, leads to a $\mathcal{P}$-comodule splitting
	\[
		\mathbb{F}_2[q_0,q_1,\dots]\cong\mathbb{F}_2\oplus\bigoplus\limits_i\mathbb{F}_2\langle q_i\rangle\oplus\bigoplus\limits_{i,j}\mathbb{F}_2\langle q_iq_j\rangle\oplus\dots
	\]
\end{prop}

This splitting reduces CESS calculations to studying monomials of fixed length $k$, their CE filtration. We now further reduce the calculations by determining which CESS elements can contribute $\tau^n$-torsion.

\begin{prop}
	Any element in the CESS with $s\leq2$ which vanishes in the abutment must support a differential.
\end{prop}
\begin{proof}
	Since all even differentials in the CESS are trivial, the first differential which can possibly be nontrivial is $d_3$, which decreases the $k$-degree by 2. Since CESS differentials raise the $(s+k)$-degree by 1, an element in total Adams filtration $s+k$ can only receive a differential from an element in total Adams filtration at most $s+k-1$. Combining this information, we find that
	\[
		d_r:E_r^{s+k-2\leq k,k+r-1,t}\to E_r^{s+k-1\leq k+1,k,t}
	\]
must be trivial, as the assumption $s\leq2$ will force the source tridegree to be identically zero. Then no element as in the proposition can receive a differential. Since such an element is killed by assumption, it must support a differential.
\end{proof}

\begin{cor}
	Any CESS element not surviving to the classical Adams $E_2$ page must support a $d_3$ differential if it is contained in $E_2^{3,3,t}$, $E_2^{4,3,t}$, or $E_2^{5,3,t}$ or must support a $d_3$ or $d_5$ differential if it is contained in $E_2^{4,4,t}$ or $E_2^{5,4,t}$.
\end{cor}
\begin{proof}
Such an element lies in the range addressed by the previous proposition, so we know that it must support a CESS differential. Since the lowest possible CE filtration is 0 and the CESS differential $d_r$ lowers this filtration by $r-1$, we find that the largest possible value of $r$ for such a differential is 3 in the first set of cases and 4 in the second set of cases.
\end{proof}

\begin{cor}
No element in CE filtration 0 or 1 supports any CESS differential. In particular, if such an element does not survive to the classical Adams $E_2$ page, then it must receive a CESS differential.
\end{cor}
\begin{proof}
Since all CESS differentials lower CE filtration by at least 2 and the lowest such filtration is 0, we see that no element in CE filtration less than 2 can support a differential.
\end{proof}

Combined, these results greatly constrain possible CESS differentials in the range $s+k\leq5$. In particular, of those elements not surviving to the Adams $E_2$ page, those in CE filtration 1 or 2 must receive $d_3$ or $d_5$ differentials, and those in CE filtration 3 or 4 must support $d_3$ differentials. As a result, out of these elements, we must have the following differential pattern, where the middle filtration is the CE filtration $k$:
\[
	\begin{split}
		E^{\ast,4,\ast} &\xrightarrow{d_3} E^{\ast,2,\ast} \\
		E^{\ast,4,\ast} &\xrightarrow{d_5} E^{\ast,0,\ast} \\
		E^{\ast,3,\ast} &\xrightarrow{d_3} E^{\ast,1,\ast} \\
		E^{\ast,2,\ast} &\xrightarrow{d_3} E^{\ast,0,\ast}
	\end{split}
\]
This pattern leads to the following conclusions, which we will use to greatly reduce the workload:

\begin{enumerate}[(1)]
	\item We need no information pertaining to $k=5$ elements
	\item The $\tau$-torsion elements in this range are images of CESS elements in CE filtration 0, 1, or 2
	\item The $\tau^2$-torsion elements in this range are images of CESS elements in CE filtration precisely 0
\end{enumerate}

In fact, the differential pattern above reduces the workload to knowing merely the following tridegrees on the CESS $E_2$ page:

\begin{center}
	\begin{tabular}{c}
		$E^{4,0,\ast}, E^{5,0,\ast}$ \\ \\
		$E^{3,1,\ast}, E^{4,1,\ast}$ \\ \\
		$E^{1,2,\ast}, E^{2,2,\ast}, E^{3,2,\ast}$ \\ \\
		$E^{0,3,\ast}, E^{1,3,\ast}$ \\ \\
		$E^{0,4,\ast}$
	\end{tabular}
\end{center}

\noindent However, it is easy to check that the only member of the last group is $q_0^4$, which must survive the CESS to yield $h_0^4\in\text{Ext}_\mathcal{A}$. As a result, all the $\tau^n$-torsion elements in the 4- and 5-lines must in fact be $\tau^1$-torsion.

In addition, all the $k=0$ elements are already determined since they do not support or receive algAHSS differentials by definition. Then they need not be considered when running the algAHSS, only when running the CESS. In total, the algAHSS workload has been reduced to the following very manageable list:

\[
	\begin{array}{c}
		E^{3,1,\ast}, E^{4,1,\ast} \\ \\
		E^{1,2,\ast}, E^{2,2,\ast} \\ \\
		E^{1,3,\ast}
	\end{array}
\]

\noindent Of course, only partial knowledge of these groups is required, and it will turn out that only about half of the listed tridegrees will be necessary to examine since they determine either sources or targets for CESS differentials.

Extending the arguments above to the 6-line, one finds that the only possibility for $\tau^2$-torsion elements on that line is that there are nontrivial differentials
\[
	E^{\ast,5,\ast}\xrightarrow{d_5}E^{\ast,1,\ast}
\]
or
\[
	E^{\ast,4,\ast}\xrightarrow{d_5}E^{\ast,0,\ast}
\]

As in the case of the 4-line, the only $k=5$ element living on the 5-line of the CESS $E_2$ page is $q_0^5$. Again, since this is the only possible element which can detect $h_0^5\in\text{Ext}_\mathcal{A}$, it must survive the CESS, so it cannot support a differential. Then we are concerned with differentials of the form $E^{\ast,4,\ast}\xrightarrow{d5}E^{\ast,0,\ast}$. We analyze such differentials in \Cref{sec:main}, where we prove \Cref{thm:6-line}.

\section{The Algebraic Atiyah-Hirzebruch Spectral Sequence}\label{sec:algAHSS}

In order to compute the CESS $E_2$ page, we follow \cite{BX} and use the algebraic Atiyah-Hirzebruch spectral sequence, herein denoted the algAHSS. This spectral sequence results from filtering the CESS $E_2$ page as follows. By the remarks in \Cref{sec:CESS}, we can rewrite the CESS $E_2$ page as
\[
	E_2^{s,k,t} = \text{Ext}_\mathcal{P}^s(\mathbb{F}_2,\mathbb{F}_2[q_0,q_1,\dots]),
\]
where we have chosen to use the classical version for ease of exposition. Again, the classical CESS $E_2$ page determines the motivic one and vice versa, so this choice is immaterial. By the splitting result described in \Cref{sec:CESS}, we can consider monomials in the $q_i$ of fixed length $k$. For each such monomial $q_{i_1}q_{i_2}\dots q_{i_k}$, we define the AH filtration of any element in
\[
	\text{Ext}_\mathcal{P}(\mathbb{F}_2,\mathbb{F}_2\langle q_{i_1}q_{i_2}\dots q_{i_k}\rangle)
\]
to be
\[
	i_1+i_2+\dots+i_k
\]
This filtration yields the algAHSS, which has the form
\[
	E_1 = \text{Ext}_\mathcal{P}\otimes\mathbb{F}_2[q_0,q_1,\dots]\Rightarrow\text{Ext}_\mathcal{P}(\mathbb{F}_2,\mathbb{F}_2[q_0,q_1,\dots])
\]

As in \cite{BX}, we use the convention that an element on the $E_1$ page can be written as $q_{i_1}q_{i_2}\dots q_{i_k}\cdot a$. The differentials in this sequence can then be expressed as
\[
	q_{i_1}q_{i_2}\dots q_{i_k}\cdot a\to q_{j_1}q_{j_2}\dots q_{j_k}\cdot b
\]
with $a\in\text{Ext}_\mathcal{P}^{s_1,t_1}(\mathbb{F}_2,\mathbb{F}_2)$ and $b\in\text{Ext}_\mathcal{P}^{s_2,t_2}(\mathbb{F}_2,\mathbb{F}_2)$. Just as with CESS differentials, all algAHSS differentials behave like Adams $d_1$ differentials in the sense that, if our differential is as above, then $s_1+1 = s_2$ and $t_1 = t_2$.

Unlike CESS differentials, for which there is no explicit formula, the algAHSS differentials can be deduced from the comodule structure $\psi$. Some special cases of these differentials are described in \cite{BX}.

\begin{prop}[\cite{BX}]\label{prop:diff}
	We have the following algAHSS differentials for $n\geq0$:
	\begin{enumerate}[(1)]
		\item $d_1(q_{n+1}\cdot a) = q_n\cdot h_na$
		\item $d_2(q_{n+1}^2\cdot a) = q_n^2\cdot h_{n+1} a$
		\item $d_2(q_{n+2}\cdot a) = q_n\cdot\langle h_n,h_{n+1},a\rangle$ if $h_{n+1}a = 0$
		\item $d_4(q_{n+2}^2\cdot a) = q_n^2\cdot\langle h_{n+1},h_{n+2},a\rangle$ if $h_{n+2}a = 0$
	\end{enumerate}
	Furthermore, algAHSS differentials satisfy a Leibniz rule.
\end{prop}

Rather than reproduce the proof here, we illustrate the general idea with the following example, which captures some of the common obstacles to running the algAHSS. Say we have already used the differentials stated above to determine that the algAHSS element $q_2q_3\cdot h_0h_2$ survives to $E_2$ and we would like to compute $d_2$ on it:
\[
	d_2(q_2q_3\cdot h_0h_2) =\text{ }?
\]
In order to evaluate this differential, there are two options in general. By \Cref{prop:diff}, algAHSS differentials satisfy a Leibniz rule. However, any decomposition of $q_2q_3\cdot h_0h_2$ will result in at least one factor which does not survive to $E_2$, so the Leibniz rule cannot be used. Instead, one must use the direct definition. To that end, one computes
\[
	\begin{split}
		\psi(q_2q_3) &= (1\otimes q_2+\xi_1^{2^2}\otimes q_1+\xi_2^{2^1}\otimes q_0)(1\otimes q_3+\xi_1^{2^3}\otimes q_2+\xi_2^{2^2}\otimes q_1+\xi_3^{2^1}\otimes q_0) \\
		&\equiv \xi_2^{2^2}\otimes q_1q_2+\xi_1^{2^2+2^3}\otimes q_1q_2+\xi_2^{2^1}\otimes q_0q_3+\dots\text{ mod }d_1
	\end{split}
\]
The comodule structure $\psi$ controls the algAHSS differentials in the following way. Recall that an algAHSS $d_r$ differential decreases the AH filtration by $r$. Since the starting filtration was
\[
	\text{AH}(q_2q_3\cdot h_0h_2) = 2+3 = 5,
\]
the $d_2$ differential must have target in AH filtration $5-2=3$. Then the target of the differential is the sum of the terms in $\psi(q_2q_3)$ with AH filtration 3. In order to translate
\[
	\xi_2^{2^2}\otimes q_1q_2+\xi_1^{2^2+2^3}\otimes q_1q_2+\xi_2^{2^1}\otimes q_0q_3
\]
back into an algAHSS element, we use the following facts. First, the $\xi_2^{2^2}$ appearing here is a witness to the nullhomotopy of $h_1h_2$ in $\text{Ext}_\mathcal{P}$. Indeed, it would normally detect the nullhomotopy of $h_2h_3$ when we work over $\mathcal{A}$, but we are working over the ``doubled-up'' algebra $\mathcal{P}$. In other words, the image of $\xi_2^{2^2}$ under the differential in the cobar complex for $\text{Ext}_\mathcal{P}$ is $[\xi_1^{2^2}|\xi_1^{2^3}]$. Then one can translate the first term in the expansion as
\[
	\xi_2^{2^2}\otimes q_1q_2 \leadsto q_1q_2\cdot\langle h_1,h_2,-\rangle
\]
Since the Ext component of $q_2q_3\cdot h_0h_2$ is $h_0h_2$, we insert that term into the Massey product to obtain
\[
q_1q_2\cdot\langle h_1,h_2,h_0h_2\rangle
\]
The Massey product appearing in the Ext component is known to be equal to $c_0$, so we obtain
\[
q_1q_2\cdot c_0
\]
We can similarly translate the final term in the expansion as
\[
q_0q_3\cdot\langle h_0,h_1,h_0h_2\rangle
\]
Now, the indeterminacy in any Massey product appearing in an algAHSS differential is killed by lower differentials by construction. After finding a single member of a Massey product, one can use that representative as the \textit{entire} Massey product. For example, the shuffling formula
\[
\langle h_0,h_1,h_0\rangle h_2\subset\langle h_0,h_1,h_0h_2\rangle,
\]
combined with the well-known product $\langle h_0,h_1,h_0\rangle = h_1^2$, yields
\[
q_0q_3\cdot\langle h_0,h_1,h_0h_2\rangle = q_0q_3\cdot h_1^2h_2 = 0,
\]
where we have used the Ext relation $h_1h_2 = 0$. For the middle term
\[
	\xi_1^{2^2+2^3}\otimes q_1q_2
\]
in the expansion, the cobar element actually witnesses a nullhomotopy yielding a matric Massey product. In this case, the matric Massey product splits into a sum of two ordinary Massey products:
\[
\langle
\begin{bmatrix}
	h_1 &h_2
\end{bmatrix}
,
\begin{bmatrix}
h_2 \\
h_1
\end{bmatrix}
,
h_0h_2
\rangle
=
\langle h_1,h_2,h_0h_2\rangle+\langle h_2,h_1,h_0h_2\rangle
\]
with no indeterminacy. The first Massey product here cancels with the first term in the expansion while the second can be evaluated, using the zero-indeterminacy argument above, as
\[
q_1q_2\cdot\langle h_2,h_1,h_0h_2\rangle = q_1q_2\cdot\langle h_2,h_1,h_2\rangle h_0 = q_1q_2\cdot h_0h_1h_3 = 0
\]
Combining all of this information,
\[
	d_2(q_2q_3\cdot h_0h_2) = 0
\]

We will now prove some results simplifying the work in evaluating algAHSS differentials.

\begin{prop}\label{prop:linear}
	All algAHSS differentials are Ext-linear. That is, if
	\[
		d_r(q_{i_1}\dots q_{i_k}\cdot a) = q_{j_1}\dots q_{j_k}\cdot b
	\]
and $c\in\text{Ext}_\mathcal{P}(\mathbb{F}_2,\mathbb{F}_2)$, then
	\[
		d_r(q_{i_1}\dots q_{i_k}\cdot ac) = q_{j_1}\dots q_{j_k}\cdot bc
	\]
\end{prop}
\begin{proof}
	This follows immediately from the Leibniz rule. Indeed, every element $a\in\text{Ext}_\mathcal{P}$ in the algAHSS is a permanent cycle since it has AH filtration 0. Then the Leibniz rule gives
	\[
		d_r(q_{i_1}\dots q_{i_k}\cdot ac) = d_r(q_{i_1}\dots q_{i_k}\cdot a)c
	\]
\end{proof}

Using this proposition, we can resolve extension problems needed for our determination of the $\mathbb{F}_2[\tau]$-module relations on the $\mathbb{C}$-motivic Adams $E_2$ page.

\begin{cor}
	Let $q_I\cdot a$ and $q_I\cdot ab$ be elements on the algAHSS $E_1$ page. If $q_I\cdot a$ survives to detect $A$ on the Adams $E_2$ page and $q_I\cdot ab$ survives to detect $B$ on the Adams $E_2$ page, then
	\[
		A\cdot\text{Sq}^0b = B
	\]
\end{cor}
\begin{proof}
	This is essentially the statement that the $\text{Ext}_\mathcal{P}$-module structure is preserved. That the algAHSS preserves this module structure follows from Ext-linearity above. For the CESS, we note that, by construction, it is a spectral sequence of $\text{Ext}_\mathcal{P}$-modules. The result follows by applying the doubling-up isomorphism $\text{Ext}_\mathcal{A}\to\text{Ext}_\mathcal{P}$.
\end{proof}

\begin{rem}
	For a specific example of this property, if one computes the $\mathbb{C}$-motivic Adams $E_2$-page in a range, as we do in \Cref{sec:low-dim}, one finds a number of $h_0$-towers. Since these extensions are by $h_0\in\text{Ext}_\mathcal{P}$, they can actually be identified with extensions by $h_1\in\text{Ext}_\mathcal{A}$ under the doubling-up isomorphism. Then one sees that these $h_0$-towers actually give the various $h_1$-diagonals, consisting of $\tau$-torsion elements, appearing in \cite{charts}. \qed
\end{rem}

\begin{prop}\label{prop:sq0}
	All algAHSS differentials are $\sq$-invariant. That is, if
	\[
		d_r(q_{i_1}\dots q_{i_k}\cdot a) = q_{j_1}\dots q_{j_k}\cdot b+\dots
	\]
	with $a,b\in\ext_\mathcal{P}$ and $r\geq i_1+\dots+i_k$, then
	\[
		d_r(q_{i_1+1}\dots q_{i_k+1}\cdot\sq a) = q_{j_1+1}\dots q_{j_k+1}\cdot\sq b+\dots
	\]
\end{prop}
\begin{proof}
	Recall that every algAHSS differential $d_r$ on an element $q_{i_1}\dots q_{i_k}\cdot a$ is determined by
	\[
		\psi(q_{i_1}\dots q_{i_k}) = \prod\limits_{n=1}^k\left(1\otimes q_{i_n}+\xi_1^{i_n}\otimes q_{i_n-1}+\xi_2^{i_n-1}\otimes q_{i_n-2}+\dots\right),
	\]
	and if $r\geq i_1+\dots+i_k$, then
	\[
		d_r(q_{i_1}\dots q_{i_k}\cdot a) = q_{j_1}\dots q_{j_k}\cdot[\xi_{\alpha_1}^{\beta_1}|\xi_{\alpha_2}^{\beta_2}|\dots|\xi_{\alpha_k}^{\beta_k}|\xi_{\gamma_1}^{\delta_1}|\dots|\xi_{\gamma_{s(a)}}^{\delta_{s(a)}}]+\dots
	\]
	where $[\xi_{\gamma_1}^{\delta_1}|\dots|\xi_{\gamma_{s(a)}}^{\delta_{s(a)}}]$ is a lift of $a\in\text{Ext}_\mathcal{P}$ to the cobar complex computing $\text{Ext}_\mathcal{P}$ and $s(a)$ denotes the Adams filtration or homological degree of $a$. Here we are treating all the $\alpha_i,\beta_i,\gamma_i,\delta_i$ as multi-indices to account for all possible cobar elements. If $r<i_1+\dots+i_k$, then all terms in $d_r$ would contain factors of $q_n$ with $n<0$, which must vanish. We elaborate on this subtlety below. Now the Ext components in the differential
	\[
		d_r(q_{i_1+1}\dots q_{i_k+1}\cdot \text{Sq}^0a) = q_{j_1+1}\dots q_{j_k+1}\cdot[\xi_{\alpha_1}^{2\beta_1}|\xi_{\alpha_2}^{2\beta_2}|\dots|\xi_{\alpha_k}^{2\beta_k}|\xi_{\gamma_1}^{2\delta_1}|\dots|\xi_{\gamma_{s(a)}}^{2\delta_{s(a)}}]+\dots
	\]
	are seen to be the result of applying $\text{Sq}^0$ to the Ext components in the first differential above. This latter differential follows from the fact that the powers $\beta_1,\dots,\beta_k$ depend linearly on $i_1,\dots,i_k$, by the expansion of $\psi(q_{i_1}\dots q_{i_k})$ above. If the cobar element
	\[
		[\xi_{\alpha_1}^{\beta_1}|\xi_{\alpha_2}^{\beta_2}|\dots|\xi_{\alpha_k}^{\beta_k}|\xi_{\gamma_1}^{\delta_1}|\dots|\xi_{\gamma_{s(a)}}^{\delta_{s(a)}}]
	\]
	detects $b\in\text{Ext}_\mathcal{P}$, then the cobar element
	\[
		[\xi_{\alpha_1}^{2\beta_1}|\xi_{\alpha_2}^{2\beta_2}|\dots|\xi_{\alpha_k}^{2\beta_k}|\xi_{\gamma_1}^{2\delta_1}|\dots|\xi_{\gamma_{s(a)}}^{2\delta_{s(a)}}]
	\]
	detects $\sq b\in\text{Ext}_\mathcal{P}$.
\end{proof}

\begin{rem}
	In the reasoning above, the condition that $r\geq i_1+\dots+i_k$ may seem artificial and not crucial to the proof of \Cref{prop:sq0}. The subtlety arises, for example, when evaluating an expression of the form
	\[
		q_{-1}q_{j_2}\dots q_{j_k}\cdot b,
	\]
	which must vanish simply because $q_{-1}$ is undefined. However, if one were to ``apply $\sq$'' by increasing the relevant $q$ indices and powers in the cobar expression for $b$, then the new expression
	\[
		q_0q_{j_2+1}\dots q_{j_k+1}\cdot \sq b
	\]
	could be nonzero. In other words, the edge case of $q_{-1}$ forces the first expression to vanish, but this expression \textit{formally} contains the same information as the second expression, up to an application of $\text{Sq}^0$. For an explicit example in which the edge case does not propagate along a $\sq$-family, note that $d_2(q_1\cdot h_0h_3^2) = 0$ for degree reasons. However,
	\[
		d_2(q_{1+1}\cdot \sq h_0h_3^2) = d_2(q_2\cdot h_1h_4^2) = q_0\cdot \langle h_0,h_1,h_1h_4^2\rangle = q_0\cdot p_0,
	\]
	which is nonzero on the algAHSS $E_2$ page.
\end{rem}

\begin{rem}
	This $\sq$-invariance allows us to evaluate algAHSS differentials in low stems, where we have extensive knowledge of the structure of the Adams $E_2$ page, and then propagate the differentials along $\sq$-families.
\end{rem}

\section{Proof of Main Results}\label{sec:main}

Using the features in \Cref{sec:algAHSS} of algAHSS differentials, we are able to use the known structure of the Adams $E_2$ page
\[
	\text{Ext}_\mathcal{A}(\mathbb{F}_2,\mathbb{F}_2),
\]
along with the ``doubling-up'' isomorphism
\[
	\text{Ext}_\mathcal{P}^{s,2t}(\mathbb{F}_2,\mathbb{F}_2)\cong\text{Ext}_\mathcal{A}^{s,t}(\mathbb{F}_2,\mathbb{F}_2),
\]
to prove our main results.

Before proving \Cref{thm:5-line}, we prove a general result about certain algAHSS cycles which will be used for the cases $k=2,3,4$. In order to state and prove this result, we first introduce some terminology, where we denote $q_{i_1}\dots q_{i_{|I|}}$ by $q_I$ for brevity.

\begin{defn}
	Say that an algAHSS $d_1$ cycle $\sum_{j\in J} q_{I_j}\cdot a_j$ is an \textit{overlap cycle} if $|J|>1$ and there is no proper subset $J'\subset J$ with $\sum_{j\in J'} q_{I_j}\cdot a_j$ a $d_1$ cycle. If a cycle consists of only a single term, then we will call it an \textit{ordinary cycle}.
\end{defn}

For an example of an overlap cycle, consider the $k=2$ element
\[
q_2q_3\cdot h_3 + q_1q_4\cdot h_1,
\]
which has
\[
	d_1(q_2q_3\cdot h_3 + q_1q_4\cdot h_1) = q_1q_3\cdot h_1h_3 + 0 + 0 + q_1q_3\cdot h_1h_3 = 0
\]
since we are working over $\mathbb{F}_2$. Then this element is a $d_1$ cycle even though neither individual term is a $d_1$ cycle. The attentive reader will notice that this overlap cycle is in fact a $d_1$ boundary and may therefore wonder whether all overlap cycles are boundaries. This is not the case in general, but it does hold for $k=2$ and Ext filtration $s\leq3$. This can be deduced from the following overlap cycle result, which says that we can bound the length of an overlap cycle, based on its CE filtration.

\begin{prop}\label{prop:overlap}
	If $\sum_{j\in J} q_{I_j}\cdot a_j$ is an overlap cycle of CE filtration $k$ and Ext filtration $s\leq3$, then its image on the algAHSS $E_2$ page can be taken to have $|J|\leq\lfloor\frac{k}{2}\rfloor$.
\end{prop}
\begin{proof}
	Say we have such an overlap cycle $\sum_{j\in J} q_{I_j}\cdot a_j$ containing terms
	\[
		\begin{tabular}{c c c}
			$q_{i_1}\dots q_{i_k}\cdot a_1$ &and &$q_{\ell_1}\dots q_{\ell_k}\cdot a_2$
		\end{tabular}
	\]

	Since neither of the terms above can be a cycle on its own, each one contributes a nonzero term when taking $d_1$. Without loss of generality, we may assume that these two nonzero $d_1$ contributions cancel each other, as they must be cancelled somehow and every term in the overlap cycle can be written in this form. By reordering the indices, we may further assume that these nonzero $d_1$ terms are
	\[
		\begin{tabular}{c c c}
			$q_{i_1-1}q_{i_2}\dots q_{i_k}\cdot h_{i_1-1}a_1 \neq 0$ &and &$q_{\ell_1-1}q_{\ell_2}\dots q_{\ell_k}\cdot h_{\ell_1-1}a_2 \neq 0$
		\end{tabular}
	\]
	Since these terms cancel by assumption and we are examining the $E_1$ page, we must have component-wise equality:
	\[
		\begin{split}
			q_{i_1-1}q_{i_2}\dots q_{i_k} &= q_{\ell_1-1}q_{\ell_2}\dots q_{\ell_k} \\
			h_{i_1-1}a_1 &= h_{\ell_1-1}a_2
		\end{split}
	\]

	One implication of the first equality is that all but two $q$ subscripts match between the two terms. If $i_1 = \ell_1$, then the remaining two subscripts also match, so $I_1 = I_2$. Since the Ext filtration is assumed to be at most 3, we have that $h_{i_1-1}a_1 = h_{i_1-1}a_2 \neq 0$ forces $a_1 = a_2$. See the remark below for more on this cancellation property. In summary, if $i_1 = \ell_1$, then the two terms are identical, so they would cancel each other even before taking $d_1$, meaning that we can disregard this case.

	If $i_1 \neq \ell_1$, then we find, again up to reindexing of subscripts, that
	\[
		\begin{split}
			q_{i_1}q_{i_2}\dots q_{i_k}\cdot a_1 &= q_{i_1}q_{i_2}q_{i_3}\dots q_{i_k}\cdot h_{i_2}a \\
			q_{\ell_1}q_{\ell_2}\dots q_{\ell_k}\cdot a_2 &= q_{i_1-1}q_{i_2+1}q_{i_3}\dots q_{i_k}\cdot h_{i_1-1}a
		\end{split}
	\]
	with the implicit assumption that $a_1,a_2$ have filtration at least 1 since otherwise we would return to the $i_1 = \ell_1$ case. Note that both Ext components are divisible by some $a\in\text{Ext}$. This is because the condition $h_{i_1-1}a_1 = h_{\ell_1-1}a_2$ forces $a_1$ to be divisible by $h_{\ell_1-1} = h_{i_2}$ and $a_2$ to be divisible by $h_{i_1-1}$. Again using the cancellation property present through Adams filtration 3, we conclude that $h_{i_1-1}h_{i_2}a_1' = h_{i_1-1}h_{i_2}a_2'$ forces $a_1' = a_2'$. Then, for each pair of terms cancelling in $d_1$ of the overlap cycle, we have characterized the elements from which these terms must originate.

	Note that in the Ext component we have a filtration-one element marking the subscript which is raised in the partner term. For example, in the first term above, we have an $h_{i_2}$ recording the fact that $i_2+1$ appears as a $q$ subscript in the second term. Similarly, there is an $h_{i_1-1}$ in the second Ext component recording the fact that $i_1-1+1 = i_1$ appears as a $q$ subscript in the first term. As a result, if our overlap cycle $\sum_{j\in J}q_{I_j}\cdot a_j$ contains a third term $q_{I_3}\cdot a_3$ overlapping with one of the two terms above, then $q_{I_3}\cdot a_3$ must take the form
	\[
		q_{i_1}q_{i_2+1}q_{i_3-1}q_{i_4}\dots q_{i_k}\cdot h_{i_3-1}a,
	\]
	again up to reindexing. We can repeat this argument until we have $\lfloor\frac{k}{2}\rfloor+1$ many terms in our overlap cycle, which may already have been the case when we had two terms above. In any case, note that every such term has an index shared between its $q$ and Ext components and each Ext component has the common factor $a$, so we find that
	\[
		d_1(q_{i_1}q_{i_2+1}q_{i_3}\dots q_{i_k}\cdot a) = q_{i_1-1}q_{i_2+1}q_{i_3}\dots q_{i_k}\cdot h_{i_1-1}a+q_{i_1}q_{i_2}q_{i_3}\dots q_{i_k}\cdot h_{i_2}a+\dots
	\]
	contains all the $\lfloor\frac{k}{2}\rfloor+1$ many terms we have just collected. Since the sum on the right-hand side is a $d_1$ boundary, it provides an $E_2$ relation under which we can replace the $\lfloor\frac{k}{2}\rfloor+1$ many terms by the remaining terms, of which there are at most $\lfloor\frac{k}{2}\rfloor$ many. This proves the claim.
\end{proof}

\begin{rem}
We take $s\leq3$ because we use the cancellation property to simplify the argument. This cancellation property does not hold in higher filtrations, as, for example, we have $h_0(h_0^3x_0) = h_0(e_0g_0) \neq 0$ in stem 37. This relation can be seen in \cite{charts}. Since the Ext components under consideration above involve elements of filtration at most 3 and two $h_i$ factors, they have filtration at most 5. This is the known range of the Adams $E_2$ page, where we can directly verify that the cancellation property holds. In higher filtrations where the cancellation property does not hold, it may still be possible to recover the proof using case-by-case analysis. However, the range studied above is all we need for our purposes.
\end{rem}

Using the overlap cycle result to reduce the workload, we can now prove \Cref{thm:5-line}, which we reproduce here for convenience:

\begin{thm*}
	As an $\mathbb{F}_2[\tau]$-module, we have
	\[
		\text{Ext}^{\leq5,\ast,\ast}_{\mathcal{A}^\text{mot}}(\mathbb{F}_2[\tau],\mathbb{F}_2[\tau])\cong F\oplus T,
	\]
	where
	\[
		F\cong\text{Ext}^{\leq5,\ast}_\mathcal{A}(\mathbb{F}_2,\mathbb{F}_2)[\tau]
	\]
	and $T$ is an $\mathbb{F}_2$-vector space generated by
	\[
		\begin{array}{c c c c}
			h_1^4\text{ in }(4,4,4) &h_1^2c_0\text{ in }(10,5,7) &h_3g\text{ in }(27,5,16) &h_1^4h_{j+5}\text{ in }(3+2^{j+5},5,4+2^{j+4})
		\end{array}
	\]
	in the stated $(t-s,s,w)$-tridegrees, where $w$ is the $\mathbb{C}$-motivic weight.
\end{thm*}

We prove this theorem by determining which $k=2$ and $k=3$ elements survive to the CESS and contribute $\tau^n$-torsion classes. By the reduction at the end of \Cref{sec:CESS}, this determines all the $\tau^1$-torsion classes, hence determines an $\mathbb{F}_2[\tau]$-basis through the 5-line.

\begin{prop}\label{prop:k=2}
The only $k=2$ elements of total filtration 3 and 4 surviving to the algAHSS $E_3$ page are of the following forms:
	\[
		\begin{tabular}{c c}
			$q_0^2\cdot h_{j+2}$ &$q_{j+1}^2\cdot h_j$
		\end{tabular}
	\]

	\[
		\begin{tabular}{c c c}
			$q_0^2\cdot h_{j+2}h_{k+2}$ &$q_{j+1}^2\cdot h_jh_k$ &$q_j^2\cdot h_{j+2}^2$ \\ \\
			$q_0q_{j+3}\cdot h_{j+1}^2$ &$q_jq_{j+3}\cdot h_{j+1}^2$ &$q_{j+2}q_{j+3}\cdot h_jh_{j+2}$
		\end{tabular}
	\]
\end{prop}
\begin{proof}

	First note that there are no overlap cycles in this range since we are considering Ext filtrations $s = 1,2$. In this range, \Cref{prop:overlap} implies that all such overlap cycles are boundaries. We first examine elements on the 2-line, which are all of the form
\[
q_jq_k
\]
By \Cref{prop:diff},
\[
	d_1(q_jq_k) = q_{j-1}q_k\cdot h_{j-1}+q_jq_{k-1}\cdot h_{k-1} = 0
\]
if and only if $k=j$. Indeed, either the terms are equal or both terms vanish because $j=0$ and $k=0$, which is already accounted for in the $k=j$ case. Then the only $d_1$ cycles on the 2-line are of the form
\[
q_j^2
\]
In fact, since no element in CE filtration $k=2$ Adams filtration less than $s+k=2$, none of these elements can be $d_1$ boundaries, so they must be nonzero on $E_2$. Moving to the 3-line, we study all elements of the form
\[
q_jq_k\cdot h_\ell
\]
By \Cref{prop:diff},
\[
	d_1(q_jq_k\cdot h_\ell) = q_{j-1}q_k\cdot h_{j-1}h_\ell+q_jq_{k-1}\cdot h_{k-1}h_\ell = 0
\]
if and only if $k=j$ or we have at least one condition from each column below:
\[
	\begin{tabular}{c | c}
		$j=0$ &$k=0$ \\ \\
		$\ell\in\{j,j-2\}$ &$\ell\in\{k,k-2\}$
	\end{tabular}
\]
Considering the four possibilities consisting of all pairs, one obtains the following $d_1$ cycles:
\[
	\begin{tabular}{c c c}
		$q_0^2\cdot h_j$ &$q_0q_j\cdot h_j$ &$q_0q_{j+2}\cdot h_j$ \\ \\
		&$q_j^2\cdot h_k$ &
	\end{tabular}
\]
The second family above consists entirely of $d_1$ boundaries, so this family is zero on $E_3$. There are members of the other families are zero on $E_3$ for the same reason, but we will eliminate these elements when we finalize the $E_3$ page below. We now run the same analysis for all elements in the 4-line, which are of the form

\[
q_jq_k\cdot h_\ell h_m
\]

\noindent Once again, such an element is a $d_1$ cycle if and only if $k=j$ or we have at least one condition from each column below:

\[
	\begin{tabular}{c | c}
		$j=0$ &$k=0$ \\ \\
		$\{\ell,m\}\cap\{j,j-2\}\neq\emptyset$ &$\{\ell,m\}\cap\{k,k-2\}\neq\emptyset$ \\ \\
		$\ell = m = j+1$ &$\ell = m = k+1$
	\end{tabular}
\]

\noindent This yields, after discarding families which consist entirely of $d_1$ boundaries, the following possibilities:

\[
	\begin{tabular}{c c c}
		$q_j^2\cdot h_kh_\ell$ &$q_0q_{j+2}\cdot h_jh_k$ &$q_0q_j\cdot h_{j+1}^2$ \\ \\
		$q_0q_{j+3}\cdot h_jh_{j+2}$ &$q_{j+2}q_{k+2}\cdot h_jh_k$ &$q_jq_{j+3}\cdot h_{j+1}^2$ \\ \\
		$q_{j+3}q_{j+4}\cdot h_jh_{j+2}$ &$q_{j+2}q_{j+3}\cdot h_jh_{j+2}$ &$q_jq_{j+3}\cdot h_jh_{j+2}$
	\end{tabular}
\]

We now evaluate $d_2$ differentials in this range. For this, one may use May spectral sequence calculations or Bob Bruner's Ext software to evaluate Massey products in low stems and then apply \Cref{prop:sq0} to evaluate $d_2$ on the entire $\sq$-family. First, we look at the 2-line elements:
\[
q_j^2
\]
Every such element is of a form described in \Cref{prop:diff}, so we compute
\[
	d_2(q_j^2) = q_{j-1}^2\cdot h_j
\]

Either $j=0$ or the target of this differential is nonzero on $E_2$. First, note that
\[
	d_2(q_0^2\cdot h_j) = 0
\]
since the $q$ index is 0 here. Since $d_1(q_0q_1) = q_0^2\cdot h_0$ and $d_2(q_1^2) = q_0^2\cdot h_1$, the only members of this family surviving to $E_3$ are of the form
\[
	q_0^2\cdot h_{j+2}
\]
for $j\geq0$. We eliminated the second family already, so we turn to the third family:
\[
	d_2(q_0q_{j+2}\cdot h_j) = q_0q_j\cdot\langle h_j,h_{j+1},h_j\rangle = q_0q_j\cdot h_{j+1}^2
\]

Note that $q_0q_j\cdot h_{j+1}^2$ is nonzero on the $E_2$ page for $j\geq0$, so all members of the $q_0q_{j+2}\cdot h_j$ family must be zero on $E_3$. We turn now to the final family. The differential
\[
	d_2(q_j^2\cdot h_k) = q_{j-1}^2\cdot h_jh_k = 0
\]
if and only if $j=0$ or $k\in\{j+1,j-1\}$. Then we see that the only $d_2$ cycles in this family are of the form
\[
	\begin{tabular}{c c c}
		$q_0^2\cdot h_{j+2}$ &$q_{j+1}^2\cdot h_j$ &$q_j^2\cdot h_{j+1}$
	\end{tabular}
\]
From above,
\[
	d_2(q_{j+1}^2) = q_j^2\cdot h_{j+1},
\]
so the final family can be written in this form if $k=j+1$, meaning it must be killed by $d_2$. The $q_0^2\cdot h_{j+2}$ family was already considered above, so we do not repeat it below in the list of $d_2$ cycles. The $q_{j+1}^2\cdot h_j$ family participates in no relations at all. In summary, the 3-line elements surviving to $E_3$ are of the forms
\[
	\begin{tabular}{c c}
		$q_0^2\cdot h_{j+2}$ &$q_{j+1}^2\cdot h_j$
	\end{tabular}
\]
for $j\geq0$, subject to no relations. This concludes the $s+k=3$ part of the proposition. For Adams filtration $s+k=4$, we use \Cref{prop:linear} to recycle some of the work done for the 3-line. In particular,
\[
	d_2(q_j^2\cdot h_kh_\ell) = q_{j-1}^2\cdot h_jh_kh_\ell
\]
Considering all the ways in which this differential can vanish, the cycles are of the forms
\[
	\begin{tabular}{c c c}
		$q_0^2\cdot h_jh_k$ &$q_{j+1}^2\cdot h_jh_k$ &$q_j^2\cdot h_{j+2}^2$ \\ \\
		&$q_{j+2}^2\cdot h_jh_{j+2}$&
	\end{tabular}
\]
The final family vanishes by $d_1(q_{j+2}q_{j+3}\cdot h_j) = q_{j+2}^2\cdot h_jh_{j+2}$. Applying the same reasoning as above, we see that the other families survive, perhaps with small modifications to the allowed range of indices:
\[
	\begin{tabular}{c c c}
		$q_0^2\cdot h_{j+1}h_{k+1}$ &$q_{j+1}^2\cdot h_jh_k$ &$q_j^2\cdot h_{j+2}^2$
	\end{tabular}
\]
These elements survive to $E_3$ for $j\geq0$, subject to no relations. We turn now to the second family:
\[
	d_2(q_0q_{j+2}\cdot h_jh_k) = q_0q_j\cdot h_{j+1}^2h_k
\]
Running the same analysis as above, we can discard the case $k=j-3$, which forces the source to vanish, and the case $k=j$, which yields a $d_1$ boundary. We then obtain the following cycles:
\[
	q_0q_{j+2}\cdot h_j^2
\]
Because of $d_1(q_1q_2\cdot h_0) = q_0q_2\cdot h_0^2$, we modify the allowed range of $j$ and see that the elements
\[
	q_0q_{j+3}\cdot h_{j+1}^2
\]
survive to $E_3$ for $j\geq0$. Note that every member of the next family is the target of a $d_2$ differential emanating from the 3-line. Then we move to the fourth family in our list:
\[
	d_2(q_0q_{j+3}\cdot h_jh_{j+2}) = q_0q_{j+1}\cdot c_j
\]
The target of this differential is nonzero on $E_2$ since
\[
	d_1(q_0q_{j+1}\cdot c_j) = q_0q_j\cdot h_jc_j = 0
\]
by known Ext relations, along with the fact that it cannot receive a $d_1$ differential. Then the entire family under consideration is zero on $E_3$. For the fifth family,
\[
	d_2(q_{j+2}q_{k+2}\cdot h_jh_k) = q_jq_{k+2}\cdot h_{j+1}^2h_k+q_{j+2}q_k\cdot h_jh_{k+1}^2
\]
by the Leibniz rule and direct examination of the comodule structure $\psi$. The target of this differential is nonzero on $E_2$, so the entire $q_{j+2}q_{k+2}\cdot h_jh_k$ family is zero on $E_3$. We turn to the next family:
\[
	d_2(q_jq_{j+3}\cdot h_{j+1}^2) = 0
\]

This differential follows from direct examination of $\psi$. Since the source is not a $d_2$ boundary, the entire family is nonzero on $E_3$. Turning to the next family, we find that
\[
	d_2(q_{j+3}q_{j+4}\cdot h_jh_{j+2}) = q_{j+1}q_{j+4}\cdot c_j,
\]
again by the Leibniz rule. The target is nonzero on $E_2$ for all $j$, so every member of the source family is zero on $E_3$. We turn to the next family:
\[
	d_2(q_{j+2}q_{j+3}\cdot h_jh_{j+2}) = 0
\]

This differential is obtained by direct inspection of $\psi(q_{j+2}q_{j+3})$ since the Leibniz rule cannot be applied here. Examining each term in this expression, the differential must vanish for degree reasons. Since no member of the source family is a $d_2$ boundary, we see that the entire family is nonzero on $E_3$. Finally, we compute
\[
	d_2(q_jq_{j+3}\cdot h_jh_{j+2}) = q_jq_{j+1}\cdot c_j
\]

Again, direction inspection of $\psi(q_{j+2}q_{j+3})$ is required here since $q_j\cdot h_jh_{j+2}$ is zero on $E_2$. However, all but one of the terms in this expression vanish for degree reasons. The remaining term can be evaluated using $\langle h_{j+1},h_{j+2},h_jh_{j+2}\rangle\ni c_j$, yielding the differential above. The target of this differential is nonzero on $E_2$ for all $j$, so the entire source family is zero on $E_3$.
\end{proof}

\begin{prop}
	Of the elements surviving to the algAHSS $E_3$ page, we have the following CESS differentials and corresponding $\tau^n$-torsion contributions:
	\[
		\begin{split}
			d_3(q_1^2\cdot h_0h_k) &\neq 0 \\
			d_3(q_2q_3\cdot h_0h_2) &\neq 0
		\end{split}
	\]
	\[
		\begin{split}
			&\Downarrow \\
			\tau h_1^4h_j &= 0 \\
			\tau h_3g &= 0
		\end{split}
	\]
\end{prop}
\begin{proof}
	Note that the $q_0^2\cdot h_{j+2}h_{k+2}$ family cannot support any algAHSS differentials, and it also cannot receive differentials from the 3-line. This is because the only two 3-line families surviving to $E_3$ are
\[
	\begin{tabular}{c c}
		$q_0^2\cdot h_{j+1}$ &$q_{j+1}^2\cdot h_j$
	\end{tabular}
\]

The first family here cannot support any differentials, and the second family here can be shown to support a $d_4$ differential hitting $q_j^2\cdot h_{j+2}^2$. This fact follows from knowledge of the Massey product $\langle h_j,h_{j+1},h_j\rangle$. As a result, the element
\[
	q_0^2\cdot h_{j+2}h_{k+2}
\]
	must be nonzero on the CESS $E_2$ page. In fact, it must detect $h_0^2h_{j+3}h_{k+3}$ on the Adams $E_2$ page. Indeed, no $k=1$ or $k=3$ element can detect $h_0^2h_{j+3}h_{k+3}$ for degree reasons. Furthermore, no $k=4$ element can detect $h_0^2h_{j+3}h_{k+3}$ since the only $k=4$ element surviving to the 4-line of the Adams $E_2$ page is $q_0^4$.

	Finally, no $k=0$ element can detect $h_0^2h_{j+3}h_{k+3}$ for degree reasons. Indeed, the only products of $h_j$ elements which could survive to detect this element would need to be of the form
	\[
		\begin{array}{c c}
			h_0h_{j+1}^2h_{k+2} &h_0h_{j+2}h_{k+1}^2
		\end{array}
	\]
	However, if the first element were to survive to the Adams $E_2$ page, then it would detect an $h_1$-divisible element, by \Cref{prop:linear} and the doubling-up isomorphism. This would contradict the fact that $h_0^2h_{j+3}h_{k+3}$ is not $h_1$-divisible, by \cite{Chen}. A similar argument rules out the second element listed above.

	Examining the 3- and 4-line elements listed in \cite{Chen}, the only other possibilities are
	\[
		\begin{array}{c c c c c}
			e_0 &p_0 &p_0' &D_3(0) &h_\ell c_0
		\end{array}
	\]
	for some suitable $\ell\geq0$. In every case, \Cref{prop:linear} and the doubling-up isomorphism give a contradiction, by the relations listed in \cite{Chen}. Then we can rule out all $k=0$ contributions. We have thus shown that $q_0^2\cdot h_{j+2}h_{k+2}$ must survive to the Adams $E_2$ page and therefore cannot contribute any $\tau^n$-torsion.

	We now examine the following family:
\[
	q_{j+1}^2\cdot h_jh_k
\]
As just remarked,
\[
	d_4(q_{j+1}^2\cdot h_j) =
	\begin{cases}
		0, &j=0 \\
		q_j^2\cdot h_{j+2}^2, &j\geq0
	\end{cases}
\]
By Ext-linearity of algAHSS differentials, this forces
\[
	d_4(q_{j+1}^2\cdot h_jh_k) =
	\begin{cases}
		0, &j=0 \\
		q_j^2\cdot h_{j+2}^2h_k, &j\geq1
	\end{cases}
\]
In the $j\geq1$ case, the targets are nonzero on the algAHSS $E_4$ page, so the only surviving members of the family under consideration are of the form
\[
q_1^2\cdot h_0h_k
\]
	For degree reasons, these elements must in fact be nonzero on the CESS $E_2$ page. They provide (1) the $\tau$-torsion in bidegree $(s,t) = (5,40)$, giving the relation $\tau h_1^4h_5 = 0$, (2) the $\tau$-torsion in bidegree $(s,t) = (5,72)$, giving the relation $\tau h_1^4h_6 = 0$, and onward. In other words, this family propagates the relation $\tau h_1^4=0$ along the 5-line.

	By the remarks above, every member of the $q_j^2\cdot h_{j+2}^2$ family receives a $d_4$ differential, so the entire family is zero in the CESS. The next two families listed in the previous proposition can be shown to support $d_3$ differentials:
\[
	\begin{split}
		q_0q_{j+3}\cdot h_{j+1}^2 &\xrightarrow{d_3} q_0q_j\cdot c_j \\
		q_jq_{j+3}\cdot h_{j+1}^2 &\xrightarrow{d_3} q_j^2\cdot c_j
	\end{split}
\]
We then focus on the final family:
\[
	q_{j+2}q_{j+3}\cdot h_jh_{j+2}
\]
Note first that
\[
	t(q_{j+2}q_{j+3}\cdot h_jh_{j+2}) = -2+2^{j+1}+2^{j+5},
\]
so when $j=0$,
\[
q_2q_3\cdot h_0h_2
\]
is the only possible element providing the $\mathbb{C}$-motivic relation
\[
\tau h_3g = 0
\]
in bidegree $(s,t) = (5,32)$ so must be nonzero on the CESS $E_2$ page. However, the rest of the $q_{j+2}q_{j+3}\cdot h_jh_{j+2}$ family participates in algAHSS differentials. Indeed, the next element $q_3q_4\cdot h_1h_3$ must be zero on the CESS $E_2$ page since it does not detect an Adams element in bidegree $(s,t) = (4,66)$ and does not provide torsion in bidegree $(s,t) = (5,66)$. For degree reasons, the element $q_3q_4\cdot h_1h_3$ does not receive algAHSS differentials and can only support an algAHSS differential hitting one of the following elements:
\[
	\begin{tabular}{c c c}
		$q_3q_4\cdot h_0^2h_3$ &$q_1q_4\cdot h_0^2h_4$ &$q_0q_4\cdot h_0h_3^2$ \\ \\
		$q_3^2\cdot h_0^2h_4$ &$q_2q_3\cdot c_1$ &$q_1q_3\cdot h_3^3$
	\end{tabular}
\]
Clearly the first element cannot be the target since it has the same AH filtration as the source. Also, the element $q_2q_3\cdot c_1$ is not a $d_1$ cycle so is zero on $E_2$. The element $q_3^2\cdot h_0^2h_4$ would have to be the target of a $d_1$ differential, but this is not the case. This rules out the second, fourth, and fifth elements.

The remaining two elements are identified on $E_2$ by the differential
\[
	d_1(q_1q_4\cdot h_3^2) = q_0q_4\cdot h_0h_3^2 + q_1q_3\cdot h_3^3,
\]
so we must have the differential
\[
	d_3(q_3q_4\cdot h_1h_3) = q_0q_4\cdot h_0h_3^2
\]
Now, since $q_0q_4\cdot h_0h_3^2$ is nonzero on $E_3$,
\[
	d_1(q_0q_4\cdot h_0h_3^2) = d_2(q_0q_4\cdot h_0h_3^2) = 0,
\]
meaning $q_jq_{j+4}\cdot h_jh_{j+3}^2$ has the same property for all $j\geq0$. This follows from \Cref{prop:sq0} since the AH filtration of $q_0q_4$ is 4, which is greater than 2, the length of the longest differential just listed. The only $d_1$ differential involving $q_jq_{j+4}\cdot h_jh_{j+3}^2$ is
\[
	d_1(q_{j+1}q_{j+4}\cdot h_{j+3}^2) = q_jq_{j+4}\cdot h_jh_{j+3}^2+q_{j+1}q_{j+3}\cdot h_{j+3}^3,
\]
the $j=0$ case being the one identifying $q_0q_4\cdot h_0h_3^2$ and $q_1q_3\cdot h_3^3$. As a result, this differential propagates to give
\[
	d_3(q_{j+2}q_{j+3}\cdot h_jh_{j+2}) = q_{j-1}q_{j+3}\cdot h_{j-1}h_{j+2}^2 \neq 0
\]
for all $j\geq1$, so we have shown that $\tau h_3g = 0$ is the only relation coming from this $k=2$ family.
\[
q_2q_3\cdot h_0h_2
\]
\end{proof}
We now turn to the $k=3$ elements:

\begin{prop}\label{prop:k=3}
	The only $k=3$ elements of total filtration 4 surviving to the algAHSS $E_3$ page are of the following forms:
	\[
		\begin{tabular}{c c c}
			$q_0^3\cdot h_{j+2}$ &$q_0q_{j+2}^2\cdot h_{j+1}$ &$q_2q_3^2\cdot h_0$
		\end{tabular}
	\]
\end{prop}
\begin{proof}
	Similarly to the $k=2$ case, the following elements form a basis for the 3- and 4-lines of the $E_2$ page:
	\[
		\begin{tabular}{c c c}
			&$q_0q_j^2$ & \\ \\
			$q_0q_j^2\cdot h_{k+1}$ & &$q_{j+2}q_k^2\cdot h_j$
		\end{tabular}
	\]
	In the first 4-line family $q_0q_j^2\cdot h_{k+1}$, we have only the following $d_2$ cycles:
	\[
		\begin{tabular}{c c c}
			$q_0^3\cdot h_{j+2}$ &$q_0q_j^2\cdot h_{j+1}$ &$q_0q_{j+1}^2\cdot h_j$
		\end{tabular}
	\]
	Every $q_0q_j^2\cdot h_{j+1}$ element receives a $d_2$ differential from $q_0q_{j+1}^2$ so is zero on $E_3$. For degree reasons, no $q_0q_{j+1}^2\cdot h_j$ element receives differentials, so all are nonzero on $E_3$, as claimed. The second and final 4-line family above gives
	\[
		d_2(q_{j+2}q_k^2\cdot h_j) = q_{j+2}q_{k-1}^2\cdot h_jh_k + q_jq_k^2\cdot h_{j+1}^2,
	\]
	which vanishes only if the terms are individually zero or cancel each other. If both terms are zero, note that the second term is a $d_1$ cycle for all values of $k$. Furthermore, the only $d_1$ relation involving the second term occurs when $k=j+1$, in which case
	\[
		q_jq_k^2\cdot h_{j+1}^2 = q_jq_{j+1}^2\cdot h_{j+1}^2 = q_{j-1}q_{j+1}q_{j+2}\cdot h_{j-1}h_{j+1},
	\]
	which vanishes when $j=0$. When $k=j+1$, we also have that the first term in our $d_2$ differential vanishes:
	\[
		q_{j+2}q_{k-1}^2\cdot h_jh_k = q_j^2q_{j+2}\cdot h_jh_{j+1} = 0
	\]
	If the terms cancel each other, then this is caused by $d_1$ relations, but we have already considered this case. Then the only $d_2$ cycle in this family is $q_2q_3^2\cdot h_0$.
\end{proof}

\begin{prop}
	Of the elements surviving to the algAHSS $E_3$ page, we have the following CESS differential and corresponding $\tau^n$-torsion contribution:
	\[
		d_3(q_2q_3^2\cdot h_0) \neq 0
	\]
	\[
		\begin{split}
			&\Downarrow \\
			\tau h_1^2c_0 &= 0
		\end{split}
	\]
\end{prop}
\begin{proof}
We have already noted that
\[
q_2q_3^2\cdot h_0
\]
cannot be hit by further algAHSS differentials, and it can be checked by low-dimensional calculations that this element does not support further differentials, either. Then it must survive all the way to the CESS $E_2$ page. Again by low-dimensional calculations, sparsity forces this element to provide the known relation $\tau h_1^2c_0 = 0$, so it supports a CESS $d_3$.

Every member of the $q_0^3\cdot h_{j+2}$ family is a permanent algAHSS cycle and cannot receive any algAHSS differentials for degree reasons. Then the entire family is nonzero on the CESS $E_2$ page, where members of this family could only receive differentials from elements of CE filtration $k=5$. However, the first such elements only appear in Adams filtration $s+k=5$. Then the entire $q_0^3\cdot h_{j+2}$ family must survive to the Adams $E_2$ page. By \Cref{prop:linear}, the element $q_0^3\cdot h_{j+2}$ detects $h_0^3h_{j+3}$ for all $j\geq0$.

The final remaining family participates in the algAHSS $d_4$ differential
	\[
		d_4(q_0q_{j+2}^2\cdot h_{j+1}) = q_0q_j^2\cdot h_{j+2}^2,
	\]
which is nonzero for $j\geq0$. Then this family does not contribute at all to the Adams $E_2$ page.
\end{proof}

This concludes the proof of \Cref{thm:5-line}.

\section{$\tau^n$-Torsion on the 6-Line}

Using analysis similar to that in the proof of \Cref{thm:5-line}, we can now prove \Cref{thm:6-line}, which we reproduce here for convenience:

\begin{thm*}
	Through the 6-line, all $\mathbb{C}$-motivic Adams $E_2$ classes are either $\tau$-free or $\tau^1$-torsion.
\end{thm*}

By \Cref{thm:5-line}, the conclusion holds through the 5-line, so we now analyze possible $\tau^n$-torsion on the 6-line.

	By the differential pattern listed in \Cref{sec:CESS}, the highest power of $\tau^n$-torsion which can appear through the 6-line is $\tau^2$-torsion. Since the only $k=5$ element surviving to the CESS is $q_0^5$, which must survive to detect $h_0^5$ in the Adams $E_2$ page, the only way for $\tau^2$-torsion to appear on the motivic 6-line is to have a nontrivial CESS differential of the form
	\[
		E^{\ast,4,\ast}\longrightarrow E^{\ast,0,\ast}
	\]

Then we first determine which $k=4$ elements survive to the CESS 5-line.

\subsection{Ordinary $d_1$ Cycles}
We first study ordinary algAHSS $d_1$ cycles. The differential
	\[
		\begin{split}
			d_1(q_jq_kq_\ell q_m\cdot h_n) = \text{ }&q_{j-1}q_kq_\ell q_m\cdot h_{j-1}h_n+q_jq_{k-1}q_\ell q_m\cdot h_{k-1}h_n \\
			&+q_jq_kq_{\ell-1}q_m\cdot h_{\ell-1}h_n+q_jq_kq_\ell q_{m-1}\cdot h_{m-1}h_n
		\end{split}
	\]
	vanishes precisely if terms match in pairs or are identically zero in pairs. Following the proof of \Cref{thm:5-line}, the only members of this family which are $d_1$ cycles are of the following forms:
	\[
		\begin{tabular}{c c c}
			$q_j^2q_k^2\cdot h_\ell$ &$q_0q_jq_k^2\cdot h_j$ &$q_0q_{j+2}q_k^2\cdot h_j$
		\end{tabular}
	\]
	The second family consists entirely of $d_1$ boundaries, and the only relations imposed on the remaining elements by $d_1$ are of the following form:
	\[
		q_jq_{j+2}q_k^2\cdot h_j = q_{j+1}^2q_k^2\cdot h_{j+1}
	\]

\subsection{Overlap $d_1$ Cycles}
	We now study the overlap cycles. By \Cref{prop:overlap}, such a cycle can be written on the $E_2$ page as a sum of at most $\lfloor\frac{k}{2}\rfloor = \lfloor\frac{4}{2}\rfloor = 2$ many ordinary cycles. Then we consider overlap cycles of the following form:
	\[
		q_jq_kq_\ell q_m\cdot h_n + q_aq_bq_cq_d\cdot h_e
	\]
	Following the proof of \Cref{prop:overlap}, we find that this sum must be of the form
	\[
		q_{j-1}q_{k+1}q_\ell q_m\cdot h_{j-1} + q_jq_kq_\ell q_m\cdot h_k,
	\]
	meaning the differential
	\[
		\begin{split}
			d_1(q_{j-1}q_{k+1}q_\ell q_m\cdot h_{j-1}+q_jq_kq_\ell q_m\cdot h_k) = &\text{ }\underbracket[0.8pt]{q_{j-1}q_{k+1}q_{\ell-1}q_m\cdot h_{j-1}h_{\ell-1}}_{(1)} +\underbracket[0.8pt]{q_{j-1}q_{k+1}q_\ell q_{m-1}\cdot h_{j-1}h_{m-1}}_{(2)}\\
			&+\underbracket[0.8pt]{q_jq_kq_{\ell-1}q_m\cdot h_kh_{\ell-1}}_{(3)}+\underbracket[0.8pt]{q_jq_kq_\ell q_{m-1}\cdot h_kh_{m-1}}_{(4)}
		\end{split}
	\]
	must vanish. Terms (1)-(4) must cancel or vanish in pairs. We summarize the various possibilities in the following table:

	\[
		\begin{tabular}{c | c | c | c | c}
			case &(1) = (2) &(1) = (3) &(1) = (4) &(1) = 0 \\
			\hline
			implies &$\ell=m$ &$k=j-1$ &redundant &$\ell=0$ or $\ell\in\{k,k+2\}$
		\end{tabular}
	\]

	\[
		\begin{tabular}{c | c | c | c}
			case &(2) = 0 &(3) = 0 &(4) = 0 \\
			\hline
			implies &$m=0$ or $m\in\{k,k+2\}$ &$\ell=0$ or $\ell\in\{j+1,j-1\}$ &$m=0$ or $m\in\{j+1,j-1\}$
		\end{tabular}
	\]
	Here we have written ``redundant'' when we obtain a condition already discovered or a special case of such a condition. Taking all the various combinations, e.g., (1)=(2) and (3)=0=(4), we obtain the following overlap cycles:
	\[
		\begin{array}{c c}
			q_{j-1}q_{k+1}q_\ell^2\cdot h_{j-1} + q_jq_kq_\ell^2\cdot h_k &q_0q_{j-1}q_{j+1}q_{j+2}\cdot h_{j-1} + q_0q_jq_{j+1}^2\cdot h_{j+1} \\
		\end{array}
	\]
	The first cycle is a $d_1$ boundary so does not survive to $E_2$. The second is not a boundary but participates in the $d_1$ relation
	\[
		d_1(q_0q_jq_{j+1}q_{j+2}) = q_0q_{j-1}q_{j+1}q_{j+2}\cdot h_{j-1}+q_0q_j^2q_{j+2}\cdot h_j+q_0q_jq_{j+1}^2\cdot h_{j+1},
	\]
	so it is already accounted for in our list of $E_2$ basis elements from the previous subsection.

\subsection{Finding $d_2$ Cycles}
	Using the presentation of $E_2$ just obtained, we determine which elements are $d_2$ cycles. Again following the proof of \Cref{thm:5-line}, the only $d_2$ cycles are of the following forms:
	\[
		\begin{tabular}{c c c c}
			$q_j^4\cdot h_k$ &$q_0^2q_j^2\cdot h_{j+1}$ &$q_0^2q_{j+1}^2\cdot h_j$ &$q_j^2q_{j+2}^2\cdot h_{j+1}$
		\end{tabular}
	\]
	The second family consists entirely of $d_2$ boundaries, so we may discard it. Again one finds that there is only one form of relations imposed by $d_2$:
	\[
		q_j^2q_{j+2}^2\cdot h_{j+1} = q_{j+1}^4\cdot h_{j+2}
	\]
	Then we may discard the final family of $d_2$ cycles, remembering that it has been incorporated into the first family under these relations. As a result, the first and third families above constitute a basis for the $E_3$ page:
	\[
		\begin{array}{c c}
			q_j^4\cdot h_k &q_0^2q_{j+1}^2\cdot h_j
		\end{array}
	\]

\subsection{Narrowing Down Surviving Families}
	The $k=4$ elements contributing $\tau^n$-torsion on the $\mathbb{C}$-motivic Adams 6-line must be of a form described above. In order to determine which of these elements do indeed make a contribution, we narrow down the families even further by analyzing higher algAHSS differentials.

	Looking first at the $q_j^4\cdot h_k$ family, we compute
	\[
		d_4(q_j^4\cdot h_k) = q_{j-1}^4\cdot h_{j+1}h_k,
	\]
	which vanishes if $j=0$ or $k\in\{j,j+2\}$ or the target is hit on an earlier page. Then we claim that the $d_4$ cycles in this family are given by
	\[
		\begin{array}{c c c c}
			q_1^4\cdot h_0 &q_0^4\cdot h_{j+2} &q_{j+1}^4\cdot h_{j+1} &q_j^4\cdot h_{j+2},
		\end{array}
	\]
	Indeed, one can check that, outside of the cases $j=0$ and $k\in\{j,j+2\}$, the only case in which the target $q_{j-1}^4\cdot h_{j+1}h_k$ of the differential can possibly vanish is when $k=j-1$, giving
	\[
		\begin{split}
			d_4(q_j^4\cdot h_k) &= q_{j-1}^4\cdot h_{j+1}h_k \\
			&= q_{j-1}^4\cdot h_{j-1}h_{j+1} \\
			&= q_{j-2}q_{j-1}^2q_j\cdot h_{j-2}h_{j+1},
		\end{split}
	\]
	where in the final line we have used a relation coming from $d_1$. One can check that this is the only $d_r$ relation in which the target participates for $r\leq3$, so the target vanishes on $E_4$ precisely if $j-2<0$. We have already accounted for the $j=0$ case, so when $j=1$, we set $k=j-1=0$ to find that
	\[
		d_4(q_1^4\cdot h_0) = 0,
	\]
	as claimed. Then we have justified our list of $d_4$ cycles for the $q_j^4\cdot h_k$ family.

	Clearly the element $q_1^4\cdot h_0$ cannot be the target of a differential. Looking at the AH filtration of $q_1^4\cdot h_0$, we see that this element is a permanent cycle, so we have shown that it must survive to the CESS.

	Next, the $q_0^4\cdot h_{j+2}$ family consists entirely of permanent cycles since every member has AH filtration 0. One can check that
	\[
		d_4(q_1^4) = q_0^4\cdot h_2,
	\]
	with $q_1^4$ clearly being nonzero on $E_4$. Then the $j=0$ member of the $q_0^4\cdot h_{j+2}$ family does not survive to the CESS. However, it is easy to see that this is the only member which dies before reaching the CESS. Then $q_0^4\cdot h_{j+3}$ survives to the CESS for $j\geq0$.

	For the $q_{j+1}^4\cdot h_{j+1}$ family, we claim that all but the $j=0$ member of this family must die on $E_8$. To see this, first note that, expanding the coproduct $\psi(q_{j+1}^4)$, all differentials $d_r$ on this family will be trivial unless $r$ is a multiple of 4. As with $q_1^4\cdot h_0$, the element $q_1^4\cdot h_1$ cannot support differentials beyond $d_4$, which we have just shown to be trivial. Then this element is a permanent cycle. For $j\geq1$, we have that
	\[
		d_8(q_{j+1}^4\cdot h_{j+1}) = q_{j-1}^4\cdot\langle h_{j+1},h_{j+2},h_{j+1}\rangle = q_{j-1}^4\cdot h_{j+2}^2 \neq 0\in E_8
	\]
	Indeed, because the only subscripts appearing in the Ext component of this target are three higher than those appearing in the $q$ component, no differential $d_r$ can hit this target for $r<8$. This can be seen by expanding the coproduct $\psi(q_j^n)$ and noting that this expansion is controlled by the 2-adic valuation of $n$. Then the target $q_{j-1}^4\cdot h_{j+2}^2$ is not hit by previous differentials. Furthermore, it can be checked directly that this target is a $d_r$ cycle for $r<8$, so it must survive to $E_8$ for all $j\geq1$. Then we have shown that the $d_8$ differential above is nonzero for $j\geq1$, so indeed the only surviving member of the $q_{j+1}^4\cdot h_{j+1}$ family is
	\[
		q_1^4\cdot h_1,
	\]
	as desired.

	Now, looking at the $q_j^4\cdot h_{j+2}$ family, note that
	\[
		d_4(q_{j+1}^4) = q_j^4\cdot h_{j+2}
	\]
	for all $j\geq0$, generalizing the differential hitting $q_0^4\cdot h_2$ above. Certainly the source $q_{j+1}^4$ of this differential survives to $E_4$ for all $j\geq0$, so we see that the entire $q_j^4\cdot h_{j+2}$ family consists of $d_4$ boundaries.

	For the remaining family $q_0^2q_{j+1}^2\cdot h_j$ found in the previous subsection, we may expand the coproduct $\psi(q_0^2q_{j+1}^2)$ to find that this family consists of $d_3$ cycles, so we compute
	\[
		d_4(q_0^2q_{j+1}^2\cdot h_j) = q_0^2q_{j-1}^2\cdot\langle h_j,h_{j+1},h_j\rangle = q_0^2q_{j-1}^2\cdot h_{j+1}^2
	\]
	Again by examining the coproduct structure, we see that the target can support a nontrivial $d_r$ differential only if $r$ is even, but
	\[
		d_2(q_0^2q_{j-1}^2\cdot h_{j+1}^2) = q_0^2q_{j-2}^2\cdot h_{j-1}h_{j+1}^2 = 0
	\]
	since we have the Ext relation $h_nh_{n+2}^2 = 0$. Then this target is a $d_r$ cycle for $r\leq3$. Using knowledge of the bases for the $E_r$ pages already obtained, one can also check directly that $q_0^2q_{j-1}^2\cdot h_{j+1}^2$ is not killed by any differential before the $E_4$ page. Then the differential obtained above kills $q_0^2q_{j+1}^2\cdot h_j$ for $j\geq0$. Note that, since the Ext component is given by $h_j$, which is not defined for $j<0$, this means the entire family must die.

	This completes the analysis of the $k=4$ elements surviving to the CESS 5-line. In summary, the following elements are the only ones that survive:
	\[
		\begin{tabular}{c c c}
			$q_0^4\cdot h_{j+3}$, $j\geq0$ &$q_1^4\cdot h_0$ &$q_1^4\cdot h_1$
		\end{tabular}
	\]

	Low-dimensional calculations, plotted in \Cref{fig:0-13}, show that $q_1^4\cdot h_0$ must in fact survive the CESS to detect $P^1h_1$ in the classical Adams $E_2$ page. From the same calculations, we find that the element $q_1^4\cdot h_1$ must survive the CESS to detect $P^1h_2$ in the classical Adams $E_2$ page. See the remarks in \Cref{sec:low-dim} for more on this pattern.

\subsection{The $q_0^4\cdot h_{j+3}$ Family}
	We now show that, for all $j\geq0$, the element $q_0^4\cdot h_{j+3}$ must survive the CESS to detect the element $h_0^4h_{j+4}$. To do so, it suffices to show that no other algAHSS elements living in the same bidegree can detect $h_0^4h_{j+4}$.

	For $k=1,3$, this follows easily from the fact that the internal degree
	\[
		t(q_0^4\cdot h_{j+3}) = 4(2^{0+1}-1)+2^{j+4} = 4+2^{j+4}
	\]
	is even, whereas the internal degrees
	\[
		\begin{split}
			t(q_j\cdot a) &= (2^{j+1}-1)+t(a) \\
			t(q_jq_kq_\ell\cdot a) &= (2^{j+1}-1)+(2^{k+1}-1)+(2^{j+1}-1)+t(a)
		\end{split}
	\]
	are odd for all $a\in\text{Ext}_\mathcal{P}$. Then all that remains is to show that there are no $k=0,2$ elements living in the bidegree of interest and surviving to detect $h_0^4h_{j+4}$.

	For $k=0$, this follows by direct inspection of the classical 5-line, as computed in \cite{Chen}. Indeed, we can rule out a product of $h_j$ elements, as the only such products with the correct internal degree would have to vanish. If one then examines the 3-, 4-, and 5-lines described in \cite{Chen}, the only remaining possibilities are
	\[
		\begin{array}{c c}
			h_4d_0 &h_7D_3(1)
		\end{array}
	\]
	since these elements have internal degree of the correct form. However, if the algAHSS element $h_4d_0$ were to survive to the Adams $E_2$ page, it would have to detect an element divisible by $d_1$, by \Cref{prop:linear} and the doubling-up isomorphism. Since $h_0^4h_6$ is not divisible by $d_1$, this is a contradiction. A similar argument rules out $h_7D_3(1)$, so we have settled the $k=0$ case.

	Finally, for $k=2$, we first compute the internal degree
	\[
		t(q_jq_k\cdot a) = 2^{j+1}+2^{k+1}-2+t(a)
	\]
	of a generic $k=2$ element with $a\in\text{Ext}^3_\mathcal{P}$. By \cite{Chen}, there are two possibilities for $a$ living on the 3-line, namely that
	\[
		\begin{tabular}{c c c}
			$a = c_\ell$ &or &$a = h_\ell h_mh_n$
		\end{tabular}
	\]
	for some $\ell,m,n$. If $a = c_\ell$, then
	\[
		t(q_jq_k\cdot a) = 2^{j+1}+2^{k+1}-2+2^{\ell+1}+2^{\ell+2}+2^{\ell+4} = 4+2^{n+4}
	\]
	Dividing both sides by 2 yields
	\[
		2^j+2^k-1+2^\ell+2^{\ell+1}+2^{\ell+3} = 2+2^{n+3}
	\]
	This is only possible if one of the exponents on the left-hand side is 0. Without loss of generality we may assume that $j\leq k$, so this forces
	\[
		\begin{tabular}{c c c}
			$j=0$ &or &$\ell=0$
		\end{tabular}
	\]
	If $j=0$, then
	\[
		d_1(q_jq_k\cdot c_\ell) = q_0q_{k-1}\cdot h_{k-1}c_\ell
	\]
	If this element survives to the CESS, then this differential must vanish, meaning that $k=0$ or $\ell\in\{k-4,k-3,k-1,k\}$. If $k=0$, then our element cannot have the correct internal degree. By direct inspection of the possible values of $\ell$, we also find that our element cannot have the correct internal degree. Then we have ruled out the case $0 = j\leq k$, so we may now assume that $0<j\leq k$ and $\ell=0$. In this case our element has internal degree
	\[
		t(q_jq_k\cdot c_0) = 2^{j+1}+2^{k+1}-2+22 = 2^{j+1}+2^{k+1}+2^4+4
	\]
	In order for the first three terms on the right-hand side to combine into a single power of 2, we would need $j=3,k=4$. This element does survive to $E_3$, but on that page we find that $d_3(q_3q_4\cdot c_0)$ contains a term $q_0q_4\cdot d_0$. This element does not participate in any relations coming from the algAHSS $d_1$ since the element $d_0$ is not divisible by any $h_i$. It also does not participate in any relations coming from the algAHSS $d_2$ since the element $d_0$ cannot be expressed as a threefold Massey product of the form $\langle h_j,h_{j+1},-\rangle$. As a result, the element $q_0q_4\cdot d_0$ survives to $E_3$ and cannot be identified with any other element in $E_3$. Since this is the only term in $d_3(q_3q_4\cdot c_0)$ with first component $q_0q_4$, we then conclude that
	\[
		d_3(q_3q_4\cdot c_0) \neq 0,
	\]
so $q_3q_4\cdot c_0$ does not survive to the CESS. Then we have ruled out any elements of the form $q_jq_k\cdot c_\ell$.

We now turn to elements of the form $q_jq_k\cdot h_\ell h_mh_n$, where we are again using \Cref{prop:overlap} to deduce that we need only study ordinary $d_1$ cycles. Rather than determine cycles of this form, which is substantially more difficult than the previous case due to the extra parameters, we will first narrow the list of options by a degree argument. This argument is purely combinatorial but rather tedious. Recall from above that
\[
	t(q_0^4\cdot h_{a+3}) = 4+2^{a+4},
\]
so we are looking for elements of total filtration 5 with this internal degree. Then we compute
\[
	t(q_jq_k\cdot h_\ell h_mh_n) = 2^{j+1}+2^{k+1}+2^{\ell+1}+2^{m+1}+2^{n+1}-2 = 4+2^{a+4},
\]
which gives
\[
	2^j+2^k+2^\ell+2^m+2^n = 3+2^{a+3} \tag{$*$}
\]
since every term was manifestly divisible by 2. Then an odd number of the subscripts on the left-hand side must be 0. This yields the following cases, where we have assumed $j\leq k$ and $\ell\leq m\leq n$ without loss of generality:
\[
	\begin{tabular}{c c c}
		\underline{\textbf{Case 1}} $q_0q_k\cdot h_\ell h_mh_n$ with $k,\ell>0$ &&\underline{\textbf{Case 2}} $q_jq_k\cdot h_0h_mh_n$ with $j,m>0$ \\ \\
		\underline{\textbf{Case 3}} $q_0^2\cdot h_0h_mh_n$ with $m>0$ &&\underline{\textbf{Case 4}} $q_jq_k\cdot h_0^3$ with $j>0$ \\ \\
		&\underline{\textbf{Case 5}} $q_0^2\cdot h_0^3$ &
	\end{tabular}
\]
Case 5 is impossible for degree reasons. Also, case 3 can be seen to be a $d_1$ boundary:
\[
	d_1(q_0q_1\cdot h_mh_n) = q_0^2\cdot h_0h_mh_n
\]
We now proceed depth-first:\newline

\noindent \underline{\textbf{Case 1:}} Starting with ($*$), we compute
\[
	2^k+2^\ell+2^m+2^n = 2(1+2^{a+2}),
\]
which has all terms divisible by 2, giving
\[
	2^{k-1}+2^{\ell-1}+2^{m-1}+2^{n-1} = 1+2^{a+2}, \tag{$\dag$}
\]
which is possible only if an odd number of the exponents on the left-hand side are 0. This yields the following cases:
\[
	\begin{tabular}{c c}
		\underline{\textbf{Case 1(a)}} $q_0q_1\cdot h_\ell h_mh_n$ with $\ell>1$ &\underline{\textbf{Case 1(b)}} $q_0q_k\cdot h_1h_mh_n$ with $k,m>1$
	\end{tabular}
\]
\vspace{0.05cm}
\[
	\underline{\textbf{Case 1(c)}}\text{ } q_0q_1\cdot h_1^2h_n\text{ with }n>1
\]
Once again, we note that the family in 1(c) consists entirely of $d_1$ boundaries, so we discard this option.\newline

\noindent \underline{\textbf{Case 1(a):}} Starting with ($\dag$), we compute
\[
	2^{\ell-1}+2^{m-1}+2^{n-1} = 2^{a+2},
\]
which is possible only if $\ell = m = n-1 = a+1$. This gives the family $q_0q_1\cdot h_{a+1}^2h_{a+2}$, which is 0 by the Ext relation $h_jh_{j+1} = 0$.\newline

\noindent\underline{\textbf{Case 1(b):}} Starting with ($\dag$), we compute
\[
	2^{k-1}+2^{m-1}+2^{n-1} = 2^{a+2},
\]
which is possible only if $k = m = n-1 = a+1$. This gives the family $q_0q_{a+1}\cdot h_1h_{a+1}h_{a+2}$, which is 0 by the same Ext relation. This settles case 1(b) and therefore settles all of case 1.\newline

\noindent\underline{\textbf{Case 2:}} Starting with ($*$), we compute
\[
	2^j+2^k+2^m+2^n = 2(1+2^{a+2}),
\]
which has all terms divisible by 2, giving
\[
	2^{j-1}+2^{k-1}+2^{m-1}+2^{n-1} = 1+2^{a+2},
\]
which is possible only if an odd number of the exponents on the left-hand side are 0. This yields the following cases:
\[
	\begin{tabular}{c c}
	\underline{\textbf{Case 2(a)}} $j=1$ and $k,m,n>1$ &\underline{\textbf{Case 2(b)}} $j=k=m=1$ and $n>1$
	\end{tabular}
\]
\vspace{0.05cm}
\[
	\underline{\textbf{Case 2(c)}}\text{ }j=m=n=1\text{ and }k>1
\]

We immediately rule out 2(b) and 2(c) since both cases lead to an Ext component in $q_jq_k\cdot h_0h_mh_n$ which vanishes by classical relations. Then we focus on 2(a). In order for
\[
	q_1q_k\cdot h_0h_mh_n
\]
to survive to the Adams $E_2$ page, it must first be a $d_1$ cycle, meaning
\[
	d_1(q_1q_k\cdot h_0h_mh_n) = q_0q_k\cdot h_0^2h_mh_n+q_1q_{k-1}\cdot h_0h_{k-1}h_mh_n = 0,
\]
again with the constraint that $k,m,n>1$ since we are working in the case 2(a). Both terms must individually vanish since there is no room for the terms to cancel. Focusing on the first term,
\[
	q_0q_k\cdot h_0^2h_mh_n = 0
\]
precisely if $m=n=2$ or $m=n=3$ since we assumed $m,n>1$. However, if $m=n=2$, then the source element $q_1q_k\cdot h_0h_mh_n$ itself vanishes, so we discard this possibility. Applying $m=n=3$ to the second term, we obtain
\[
	q_1q_{k-1}\cdot h_0h_{k-1}h_3^2 = 0
\]
precisely if $k\in\{2,4,5\}$. This follows from Chen's relations through the Adams 5-line in \cite{Chen}. All of this analysis yields the following candidates:
\[
	\begin{array}{c c c}
		q_1q_2\cdot h_0h_3^2 &q_1q_4\cdot h_0h_3^2 &q_1q_5\cdot h_0h_3^2
	\end{array}
\]

Only the internal degree of the second element is of the correct form, so it is the only one we consider. Suppose, then, that $q_1q_4\cdot h_0h_3^2$ survives to the Adams $E_2$ page to detect $h_0^4h_6$. Then $q_0^4\cdot h_5$ cannot survive to the Adams $E_2$ page since there is only a single $\mathbb{C}$-motivic generator in that bidegree, which we just assumed to be detected by $q_1q_4\cdot h_0h_3^2$.

This would force $q_0^4\cdot h_5$, which survives to the CESS by the reasoning above, to participate in a CESS differential. However, if one inspects \cite{charts} in this region, one finds that there is no $\tau^n$-torsion appearing in the appropriate bidegree for this to be the case. Then $q_0^4\cdot h_5$ must survive to the Adams $E_2$ page and therefore detect the sole $\mathbb{F}_2$-generator $h_0^4h_6$ in that bidegree.

\section{Low-Dimensional Calculations}\label{sec:low-dim}

In this section we provide charts of the first 29 stems of the $\mathbb{C}$-motivic Adams $E_2$ page. We also discuss some of the calculations behind these charts and various patterns visible in the algAHSS representatives which detect the resulting Adams elements. Calculations have been made into stem 30, but the algAHSS differentials around Adams filtration 9 require far more knowledge of Massey products on the Adams $E_2$ page. We comment on stem 30 but do not show it in our charts because our calculations are incomplete at the moment.

\begin{figure}
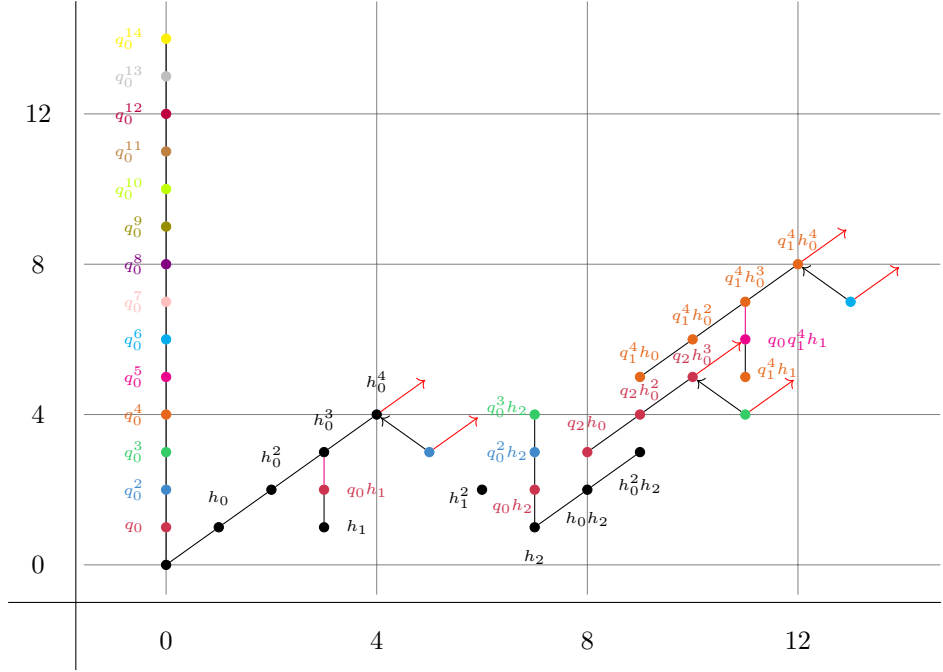

	\begin{sseqpage}[Adams grading,classes = {fill,font = \tiny},grid = go,xscale = 0.7,yscale = 0.5,xrange = {-1}{14},x tick step = 4,y tick step = 4,grid step = 4]
		\class(0,0)
		\class[pastelred,"q_0" left](0,1)
		\class[pastelblue,"q_0^2" left](0,2)
		\class[pastelgreen,"q_0^3" left](0,3)
		\class[pastelorange,"q_0^4" left](0,4)
		\class[magenta,"q_0^5" left](0,5)
		\class[cyan,"q_0^6" left](0,6)
		\class[pink,"q_0^7" left](0,7)
		\class[violet,"q_0^8" left](0,8)
		\class[olive,"q_0^9" left](0,9)
		\class[lime,"q_0^{10}" left](0,10)
		\class[brown,"q_0^{11}" left](0,11)
		\class[purple,"q_0^{12}" left](0,12)
		\class[lightgray,"q_0^{13}" left](0,13)
		\class[yellow,"q_0^{14}" left](0,14)
		\foreach \n in {0,...,13}{
			\structline(0,\n)(0,\n+1)
		}

		\class["h_0" above](1,1)
		\structline(0,0)(1,1)
		\class["h_0^2" above](2,2)
		\structline(1,1)(2,2)
		\class["h_0^3" above](3,3)
		\structline(2,2)(3,3)
		\class["h_0^4" above](4,4)
		\structline(3,3)(4,4)
		\class[draw = none,fill = none](5,5)
		\structline[->,red](4,4)(5,5)
		\class["h_1" right](3,1)
		\class[pastelred,"q_0h_1" right](3,2)
		\structline(3,1)(3,2)
		\structline[magenta](3,2)(3,3)
		\class[pastelblue](5,3)
		\class[draw = none,fill = none](6,4)
		\structline[->,red](5,3)(6,4)
		\d1(5,3)

		\class["h_1^2" {left,xshift=4pt,yshift=-3pt}](6,2)
		\class["h_2" below](7,1)
		\class[pastelred,"q_0h_2" {left,xshift=8pt,yshift=-6pt}](7,2)
		\structline(7,1)(7,2)
		\class[pastelblue,"q_0^2h_2" {left,xshift=6pt}](7,3)
		\structline(7,2)(7,3)
		\class[pastelgreen,"q_0^3h_2" {left,xshift=6pt,yshift=3pt}](7,4)
		\structline(7,3)(7,4)

		\class["h_0h_2" below](8,2)
		\structline(7,1)(8,2)
		\class["h_0^2h_2" below](9,3)
		\structline(8,2)(9,3)

		\class[pastelred,"q_2h_0" above](8,3)
		\class[pastelred,"q_2h_0^2" {above,yshift=-3pt}](9,4)
		\structline(8,3)(9,4)
		\class[pastelred,"q_2h_0^3" {above,yshift=-4pt}](10,5)
		\structline(9,4)(10,5)

		\class[pastelorange,"q_1^4h_0" {above,yshift=-3pt}](9,5)
		\class[pastelorange,"q_1^4h_0^2" {above,yshift=-3pt}](10,6)
		\structline(9,5)(10,6)
		\class[pastelorange,"q_1^4h_0^3" {above,yshift=-3pt}](11,7)
		\structline(10,6)(11,7)
		\class[pastelorange,"q_1^4h_1" {right,xshift=-4pt,yshift=3pt}](11,5)
		\class[magenta,"q_0q_1^4h_1" right](11,6)
		\structline(11,5)(11,6)
		\structline[magenta](11,6)(11,7)
		\structline[->,red](10,5)(11,6)
		\class[pastelorange,"q_1^4h_0^4" {above,yshift=-3pt}](12,8)
		\class[draw = none,fill = none](13,9)
		\structline[->,red](12,8)(13,9)
		\structline(11,7)(12,8)

		\class[pastelgreen](11,4)
		\class[draw = none,fill = none](12,5)
		\structline[->,red](11,4)(12,5)
		\d1(11,4)

		\class[cyan](13,7)
		\class[draw = none,fill = none](14,8)
		\structline[->,red](13,7)(14,8)
		\d1(13,7)
	\end{sseqpage}
	\caption{Stems 0 Through 13}
	\label{fig:0-13}
\end{figure}

We plot our calculations using three separate charts so that more of the structure can be seen. All $h_0$- and $h_1$-extensions present on the $\mathbb{C}$-motivic Adams $E_2$ page are shown in our charts. In \Cref{fig:0-13} we label most of the classes to demonstrate how to translate between the algAHSS names for these classes and their $E_2$ names. In the latter two charts we label only the $h_0$- and $h_1$-generators, as the structlines can be used to deduce the names of the other classes using the pattern in \Cref{fig:0-13}.

\begin{figure}
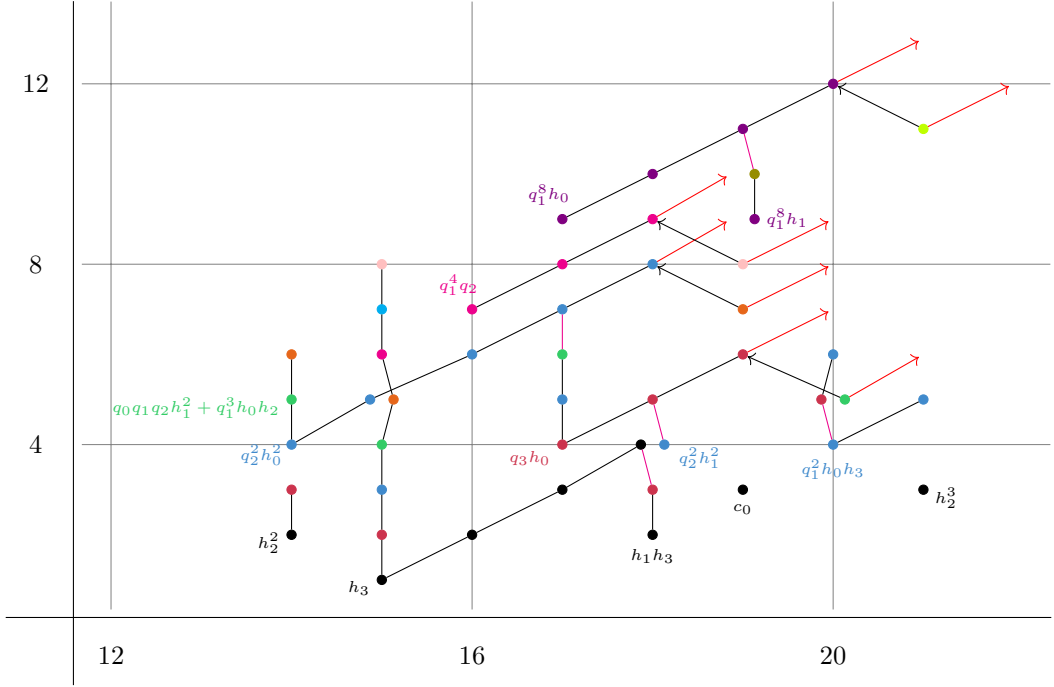

	\begin{sseqpage}[Adams grading,classes = {fill,font = \tiny},grid = go,xscale = 1.2,yscale = 0.6,xrange = {12}{22},x tick step = 4,y tick step = 4,grid step = 4]

		\class["h_2^2" {left,xshift=4pt,yshift=-3pt}](14,2)
		\class["h_2^3" {right,xshift=-4pt,yshift=-3pt}](21,3)

		\class["h_3" {left,xshift=4pt,yshift=-3pt}](15,1)
		\class(16,2)
		\class(17,3)
		\class(18,4)

		\class["h_1h_3" {below,yshift=3pt}](18,2)

		\class["c_0" {below,yshift=3pt}](19,3)

		\class[pastelred](15,2)
		
		\class[pastelred](18,3)
		\class[pastelred](14,3)

		\class[pastelred](20,5)

		\class[pastelred, "q_3h_0" {left,xshift=4pt,yshift=-5pt}](17,4)
		\class[pastelred](18,5)
		\class[pastelred](19,6)
		\class[draw=none,fill=none](20,7)
		\structline[->,red](19,6)(20,7)

		\class[pastelblue](15,3)

		\class[pastelblue](20,6)

		\class[pastelblue, "q_1^2h_0h_3" {below,yshift=3pt}](20,4)

		\class[pastelblue](21,5)

		\class[pastelblue](17,5)

		\class[pastelblue, "q_2^2h_0^2" {left,xshift=5pt,yshift=-4pt}](14,4)
		\class[pastelblue, "q_2^2h_1^2" {right,xshift=-3pt,yshift=-5pt}](18,4)

		\class[pastelblue](15,5)
		\class[pastelblue](16,6)
		\class[pastelblue](17,7)
		\class[pastelblue](18,8)
		\class[draw = none,fill = none](19,9)
		\structline[->,red](18,8)(19,9)

		\class[pastelgreen](15,4)

		\class[pastelgreen](17,6)

		\class[pastelgreen](20,5)
		\class[draw=none,fill=none](21,6)
		\structline[->,red](20,5,2)(21,6)

		\class[pastelorange](15,5)

		\class[pastelorange](19,7)
		\class[draw = none,fill = none](20,8)
		\structline[->,red](19,7)(20,8)

		\class[magenta](15,6)

		\class[magenta, "q_1^4q_2" {above,xshift=-5pt,yshift=-3pt}](16,7)
		\class[magenta](17,8)
		\class[magenta](18,9)
		\class[draw = none,fill = none](19,10)
		\structline[->,red](18,9)(19,10)

		\class[cyan](15,7)
		\class[pink](15,8)

		\class[violet, "q_1^8h_0" {above,xshift=-5pt,yshift=-3pt}](17,9)
		\class[violet](18,10)
		\class[violet](19,11)
		\class[violet, "q_1^8h_1" {right,xshift=-4pt}](19,9)
		\class[olive](19,10)
		\class[violet](20,12)
		\class[draw = none,fill = none](21,13)
		\structline[->,red](20,12)(21,13)

		%still need to fix?
		\class[pastelgreen, "q_0q_1q_2 h_1^2 + q_1^3 h_0h_2" {left,xshift=4pt,yshift=-3pt}](14,5)
		\class[pastelorange](14,6)

		\class[pink](19,8)
		\class[draw = none,fill = none](20,9)
		\structline[->,red](19,8)(20,9)
		\d1(19,8)

		\class[lime](21,11)
		\class[draw = none,fill = none](22,12)
		\structline[->,red](21,11)(22,12)
		\d1(21,11)

		\d1(20,5,2)

		\d1(19,7)

		%STRUCTLINES
		\structline(14,2)(14,3)

		\structline(14,4)(14,5)
		\structline(14,5)(14,6)

		\foreach \n in {1,2,3,6,7}{
			\structline(15,\n)(15,\n+1)
		}
		\structline(15,4)(15,5,2)
		\structline(15,5,2)(15,6)

		\foreach \n in {14,...,17}{
			\structline(\n,\n-10)(\n+1,\n-9)
		}

		\structline(15,1)(16,2)
		\structline(16,2)(17,3)
		\structline(17,3)(18,4)

		\structline(16,7)(17,8)
		\structline(17,8)(18,9)

		\structline(17,9)(18,10)
		\structline(18,10)(19,11)
		\structline(19,9,2)(19,10,2)
		\structline[magenta](19,10,2)(19,11)
		\structline(19,11)(20,12)

		\structline(17,4)(17,5)
		\structline(17,5)(17,6)
		\structline[magenta](17,6)(17,7)

		\structline(18,2)(18,3)
		\structline[magenta](18,3)(18,4)

		\structline[magenta](18,4,2)(18,5)
		\structline(17,4)(18,5)
		\structline(18,5)(19,6)

		\structline(20,4)(21,5)
		\structline[magenta](20,4)(20,5)
		\structline(20,5)(20,6)
	\end{sseqpage}
	\caption{Stems 14 Through 21}
	\label{fig:14-21}
\end{figure}

For example, the red class in bidegree $(t-s,s)=(14,3)$ in \Cref{fig:14-21} is an $h_0$-extension of $h_2^2$, so by the patterns in \Cref{fig:0-13}, we reason that the red class must be detected by $q_0\cdot h_2^2$. This does not always hold, however, and we mark generators that deviate from this pattern. Looking at the $h_1$-ladder beginning in bidegree $(23,4)$ in \Cref{fig:22-29}, one sees that the first and third generator have been labelled. This is because the third generator is not $q_2\cdot h_0^4h_3$ as one would obtain by applying the pattern in \Cref{fig:0-13}. Indeed, the classical relation $h_0^4h_3 = 0$ would force this element to vanish. Instead, another algAHSS element $q_1\cdot P^1h_2$ takes its place.

To demonstrate how algAHSS elements survive to detect Adams $E_2$ elements, we have color-coded all generators. Built into \Cref{fig:0-13} is a simple color legend in stem 0, which should help in navigating the charts. For example, at filtration 0, we have a black dot, so all black dots are images of $k=0$ algAHSS elements. Similarly, at filtration 1, we have a red dot, so all red dots are images of $k=1$ algAHSS elements. Continuing in this way, one can then track CE filtration in the rest of the chart.

One note of caution: the purple elements in $k=12$ look fairly similar to the red elements in $k=1$. These elements never appear next to each other in the charts, hence the choice of colors, but the reader should be aware that, e.g., the $v_1$-periodic pattern in stems 25 through 28 is detected by $k=12$ elements. This $v_1$-periodicity is discussed below.

In order to tidy up our charts, we have used the convention in \cite{charts} of showing the start of a $\tau$-torsion $h_1$-diagonal and then indicating that this pattern continues along the direction of a red arrow. Since our charts also show the elements \textit{providing} the $\tau$-torsion, we have indicated that those source elements continue along $h_1$-diagonals as well. The negative-slope black arrows then correspond to CE $d_3$ differentials, which drop CE filtration by $3-1=2$. If one were to compute a larger range of stems, higher CE differentials would also occur, but in this range there are only $d_3$ differentials.

\begin{figure}
	\begin{sseqpage}[Adams grading,classes = {fill,font = \tiny},grid = go,xscale=1.1,yscale = 0.8,xrange = {20}{30},x tick step = 4,y tick step = 4,grid step = 4]

		\class(22,4)

		\class["P^1h_1" {below,yshift=3pt}](23,5)
		\class(27,5)

		%\class(30,2)

		\class["h_1P^1h_1" {below,yshift=5pt}](26,6)
		\class(29,7)
		\class(28,6)

		%\class[pastelred](30,3)

		\class[pastelred](23,6)
		
		\class[pastelred](26,7)

		\class[pastelred, "q_1P^1h_2" {left,xshift=7pt,yshift=-5pt}](25,6)

		\class[pastelred, "q_2h_0^2h_3" {below}](23,4)

		\class[pastelred](24,5)

		%\class[pastelblue](30,4)

		\class[pastelblue, "q_1q_2h_1c_0" {left}](30,6)

		\class[pastelblue, "q_2q_3h_0h_2" {below}](28,4)
		
		\class[pastelblue](29,5)

		\class[pastelblue, "q_2q_3 h_0^5" {left,xshift=7pt,yshift=-5pt}](25,7)
		\class[draw = none,fill = none](26,8)
		\structline[->,red](25,7)(26,8)

		%\class[pastelgreen](30,5)

		%\class[pastelgreen](30,7)

		\class[pastelgreen, "q_0q_2q_3 h_1^2" {below,xshift=-5pt,yshift=3pt}](26,5)

		\class[pastelgreen, "q_1^2q_2 h_0^2h_3" {above,xshift=9pt,yshift=-3pt}](27,6)

		\class[pastelgreen, "q_1^2q_3 h_0^3" {below}](22,7)
		\class[draw = none,fill = none](23,8)
		\structline[->,red](22,7)(23,8)

		%\class[pastelorange](30,8)

		\class[pastelorange, "q_0q_1^2q_2 h_0h_3" {below,xshift=25pt,yshift=5pt}](26,6)
		\class[draw = none,fill = none](27,7)
		\structline[->,red](26,6,2)(27,7)

		\class[pastelorange, "q_0q_2^2q_3 h_0^3" {right}](29,7)

		\class[pastelorange, "q_2^4 h_0^4" {below,xshift=5pt,yshift=3pt}](28,8)
		\class[pastelorange](29,9)
		%\class[pastelorange](30,10)

		%\class[pastelorange, "q_2^4 h_1^2" {left}](30,6)

		%\class[magenta](30,9)

		\class[magenta, "q_0^2q_1^2q_3 h_0h_2" {right,yshift=7pt}](26,7)

		\class[pink, "q_0q_1^5q_2 h_1^2 + q_1^7 h_0h_2" {left,xshift=4pt,yshift=-3pt}](22,9)
		%\class[draw = none,fill = none](23,10)
		%\structline[->,red](22,9)(23,10)

		\class[violet](22,10)

		\class[magenta](29,8)

		\class[magenta, "q_1^5 h_0^3h_3" {left,xshift=5pt,yshift=-3pt}](28,9)

		\class[magenta, "q_1^4q_3 h_0^3" {left,xshift=7pt,yshift=-5pt}](25,8)
		\class[magenta](26,9)
		\class[magenta](27,10)
		\class[draw = none,fill = none](28,11)
		\structline[->,red](27,10)(28,11)

		\class[magenta, "q_1^2q_2^3 h_0" {right,xshift=-5pt}](23,6)
		\class[draw = none,fill = none](24,7)
		\structline[->,red](23,6,2)(24,7)

		\class[lime, "q_1^8q_2 h_0" {above}](24,11)
		\class[lime](25,12)
		\class[lime](26,13)
		\class[draw = none,fill = none](27,14)
		\structline[->,red](26,13)(27,14)

		\class[cyan, "q_1^4q_2^2 h_0^2" {below,xshift=-5pt,yshift=3pt}](22,8)

		\class[cyan, "q_1^4q_2^2 h1" {right,xshift=-3pt,yshift=5pt}](23,7)
		\class[pink](23,8)
		\class[cyan](23,9)
		\class[violet](23,9)
		\class[olive](23,10)
		\class[lime](23,11)
		\class[brown](23,12)

		\class[cyan](24,10)

		\class[cyan](25,9)
		\class[pink](25,10)
		\class[cyan](25,11)
		\class[purple, "q_1^{12} h_0" {above}](25,13)

		\class[cyan](26,8)
		\class[cyan](26,12)
		\class[draw = none,fill = none](27,13)
		\structline[->,red](26,12)(27,13)
		\class[violet](27,11)
		\class[draw = none,fill = none](28,12)
		\structline[->,red](27,11)(28,12)
		\class[purple](26,14)

		\class[purple, "q_1^{12} h_1" {right,xshift=-3pt,yshift=4pt}](27,13)
		\class[lightgray](27,14)
		\class[purple](27,15)
		\class[purple](28,16)
		\class[draw = none,fill = none](29,17)
		\structline[->,red](28,16)(29,17)
		\class[cyan](28,10)

		%\class[cyan](30,10)
		%\class[pink, "q_1^6q_3 h_0^3" {left}](30,11)

		%TORSION SOURCES
		%\class[pink](23,8)
		%\class[draw = none,fill = none](24,9)
		%\structline[->,red](23,8,2)(24,9)
		%\d1(23,8,2)

		\class[pink](28,9)
		\class[draw = none,fill = none](29,10)
		\structline[->,red](28,9,2)(29,10)
		\d1(28,9,2)

		\class[purple](27,12)
		\class[draw = none,fill = none](28,13)
		\structline[->,red](27,12)(28,13)
		\d1(27,12)

		\class[yellow](29,15)
		\class[draw = none,fill = none](30,16)
		\structline[->,red](29,15)(30,16)
		\d1(29,15)

		%DIFFERENTIALS
		\d1(28,4)
		\d1(29,5)
		\d1(30,6)

		\d1(26,5)
		\d1(27,6)

		\d1(26,6,2)

		\d1(23,6,2)

		\d1(27,11)

		%stems 20,21 for exotic extensions

		\class[pastelblue, "q_1^2 h_0h_3" {below,xshift=5pt}](20,4)
		\class[pastelred](20,5)
		\class[pastelblue](20,6)
		\class[pastelblue](21,5)

		\class(21,3)

		%STRUCTLINES
		\structline(22,8)(22,9)
		\structline(22,9)(22,10)
		\structline(22,8)(23,9)
		\structline(23,9)(24,10)
		\structline(24,10)(25,11)
		\structline(25,11)(26,12)

		\structline[magenta](20,4)(23,5) \structline[magenta](20,4)(20,5)
		\structline(20,5)(20,6)
		\structline(20,4)(21,5)
		\structline(23,5)(23,6)

		\structline(23,4)(24,5)
		\structline(24,5)(25,6)
		\structline(25,6)(26,7)

		\structline(23,7)(23,8,2)
		\structline(23,8,2)(23,9,2)
		\structline(23,9,2)(23,10)
		\structline(23,10)(23,11)
		\structline(23,11)(23,12)

		\structline(24,11)(25,12)
		\structline(25,12)(26,13)

		\structline(25,8)(25,9)
		\structline(25,9)(25,10)
		\structline[magenta](25,10)(25,11)
		\structline(25,8)(26,9)
		\structline(26,9)(27,10)

		\structline(26,6)(26,7)

		\structline(26,7,2)(26,8,2)
		\structline[magenta](26,8,2)(26,9)

		\structline(25,13)(26,14)
		\structline(26,14)(27,15)
		\structline(27,15)(28,16)

		\structline(27,5)(28,6)
		\structline(28,6)(29,7)

		\structline(27,13,2)(27,14,2)
		\structline[magenta](27,14,2)(27,15)

		\structline(28,8)(29,9)
		%\structline(29,9)(30,10)
		\structline(28,8)(28,9)
		\structline(28,9)(28,10)

		\structline(29,7,2)(29,8)
		\structline[magenta](29,8)(29,9)

	\end{sseqpage}
	\caption{Stems 20 Through 29}
	\label{fig:22-29}
\end{figure}

We now make a few comments on what these low-dimensional calculations reveal about the structure of the Adams $E_2$ page, in both the classical and $\mathbb{C}$-motivic settings.

\subsection{$h_0$-Divisibility}\label{subsec:h0}
The first thing to notice is that the $h_0$-divisibility of the various $h_j$ on the Adams $E_2$ page is enforced purely by combinatorics. In particular, the survival of any element of the form $q_0^k\cdot h_{k'}$ to the CESS is determined by the 2-adic valuation of $k$. Say $k$ has the dyadic expansion
\[
	k = 2^{i_1}+\dots+2^{i_n}
\]
with $i_1<\dots<i_n$. Then we have the algAHSS differentials
\[
	\begin{split}
		d_1(q_0^{k-1}q_1) &= q_0^k\cdot h_0 \\
		d_2(q_0^{k-2}q_1^2) &= q_0^k\cdot h_1 \\
		d_4(q_0^{k-4}q_1^4) &= q_0^k\cdot h_2 \\
		&\vdots \\
		d_{2^{i_n}}(q_0^{k-2^{i_n}}q_1^{2^{i_n}}) &= q_0^k\cdot h_{i_n}
	\end{split}
\]

Note that the sources of the differentials listed above cannot be algAHSS boundaries, as there are no elements which can possibly support differentials hitting them. Furthermore, by expanding the coproduct $\psi(q_0^{k-2^i}q_1^{2^i})$, one finds that every differential listed above is the first one supported by that source. As a result, we see that the first $q_0^k\cdot h_{k'}$ possibly surviving to the Adams $E_2$ page is $q_0^k\cdot h_{i_n+1}$. This bound is already present before running the CESS.

\subsection{Periodicity}
Another feature of \Cref{fig:0-13} is the periodicity mentioned in \Cref{sec:main}. Namely, the Adams periodicity operator $P^1(-)$ is detected by $q_1^4$, in the sense that, for example, $P^1h_{j+1}$ is detected by $q_1^4\cdot h_j$. Furthermore, the element detecting $P^2h_{j+1}$ is of the form $q_1^8\cdot h_j$ while the element detecting $P^3h_{j+1}$ is of the form $q_1^{12}\cdot h_j$.

In fact, when finishing the calculations for these charts, we were able to use this pattern to very quickly identify some of the remaining contributions from CE filtration greater than 6. Indeed, we found that, if an element $a\in\text{Ext}_\mathcal{P}$ was detected by an algAHSS element $q_I\cdot b$, then $P^1a$ was detected by $q_1^4q_I\cdot b$. Of course, not all elements can be found easily in this way, which is why there is still work to be done in stem 30.

This pattern persists beyond the low-dimensional calculations presented here, as pointed out in subsection 6.1 of \cite{BX}. There it is noted that the higher $v_i$ appearing in cobar calculations for the Adams-Novikov spectral sequence correspond to the $q_i$ in our calculations. This explains the $q_1^4$ pattern described above, as the Adams periodicity operator $P^1(-)$ is related to multiplication in $BP$-homology by $v_1^4$.

Although it is not plotted in our charts, the 30-stem element called $\Delta h_2^2$ in \cite{charts} is detected by $q_2^4\cdot h_1^2$. Recall that $h_1^2\in\text{Ext}_\mathcal{P}$ maps under the doubling-up isomorphism to $h_2^2\in\text{Ext}_\mathcal{A}$, so the form of $q_2^4\cdot h_1^2$ indicates that this $\Delta h_2^2$ should be some $v_2^4$-multiple of $h_2^2$. Indeed, the notation $\Delta h_2^2$ is meant to convey the fact that this element is detected by a multiple of an element $\Delta\in\pi_\ast tmf$. This element $\Delta$ in fact behaves as a $v_2$-periodicity operator, in particular, as multiplication by $v_2^4$.

It seems plausible that candidates for algAHSS elements detecting certain classes $a\in\text{Ext}_\mathcal{A}$ can be found fairly easily by using this pattern, along with low-dimensional calculations. One can then go about applying the program of \cite{BX} to study the survival of these classes in the Adams spectral sequence.

\subsection{Unrealized Periodicity}
One might object that there are classes plotted in our charts that are $q_1^4$- and $q_2^4$-multiples in the algAHSS but not $v_1^4$- or $v_2^4$-multiples on the Adams $E_2$ page. Focusing on the specific example of the class in bidegree $(23,6)$, we see that the class $i$ in \cite{charts} is detected by an algAHSS element called $q_1^4q_2^2\cdot h_1$. This suggests that $i$ is $P^1a$ for some $a\in\text{Ext}$, but this element $a$ would have to be $h_0^2h_4$ for degree reasons. A quick Massey product calculation shows that this is not the case.

In fact, the element $q_1^4q_2^2\cdot h_1$ is a $q_1^4$-multiple of $q_2^2\cdot h_1$, which supports an algAHSS differential
\[
	d_4(q_2^2\cdot h_1) = q_0^2\cdot h_2^2
\]
similar to those described in \Cref{subsec:h0}. This differential has the effect of preventing $h_0^2h_3^2$ from surviving to the Adams $E_2$ page. Then one might say that $i$ ``wants'' to belong to a $v_1$-periodic family starting in stem 15. In this direction, there is an algAHSS differential
\[
	d_2(q_1^6q_2) = q_0^2q_1^4q_2^2\cdot h_1+q_1^8\cdot h_2
\]
identifying $h_0^2i$ and this algAHSS element $q_1^8\cdot h_2$. Under the doubling-up isomorphism, the element $h_2\in\text{Ext}_\mathcal{P}$ maps to $h_3\in\text{Ext}_\mathcal{A}$, so this element should detect a member of the Massey product
\[
	P^2h_3
\]
on the Adams $E_2$ page. Indeed, if one computes this Massey product, one finds that it contains $h_0^2i$. Furthermore, there is a member of this Massey product which is in the image of $J$. One can see this in the charts of \textcite{Bruner-RognesJ}, where it is denoted by $h_3w_1^2$.

Since $i$ and $h_0i$ support Adams differentials, they do not survive to the image of $J$. Considering these algAHSS calculations, one might say that $i$ and $h_0i$ fail to be detected by $J$ because they are not images of the Hopf classes under repeated applications of the Adams periodicity operator $P$. The class $h_0^2i$, on the other hand, detects a member of $P^2h_3$ by the comments above, allowing $h_0^2i$ to be detected in the image of $J$.

\subsection{Exotic Extensions}
Now we will point out another feature of these calculations by first defining an extension on the $\mathbb{C}$-motivic Adams $E_2$ page to be \textit{exotic} if it is of the form $ab = \tau^n c$ for some $n\geq1$ and with neither of $a,b$ divisible by $\tau$. The first exotic extension appears in stem 3 and is of the form
\[
	h_0^2h_2 = \tau h_1^3,
\]
as seen in \cite{charts}. Classically, the corresponding relation is $h_0^2h_2 = h_1^3$, but in the $\mathbb{C}$-motivic setting there is a motivic weight which requires this factor of $\tau$. This exotic extension can be seen in our low-dimensional charts, specifically in \Cref{fig:0-13}, where we see a $q_0$-extension from $q_0\cdot h_1$ to $h_0^3$. Since the source of the extension lies in CE filtration 1 while the target lies in CE filtration 0, a multiple of $\tau$ is required to label a ``mismatch'' in filtrations.

One may be tempted to conjecture that the CE filtration serves as a proxy for motivic weight, but this is not quite true. To see this, note that we can view the extension above as an $h_2$-extension from $q_0^2$, which detects $h_0^2$, to $h_0^3$, which detects $h_1^3$. Since $q_0^2$ lies in filtration 2 while $h_0^3$ lies in filtration 0, we see that factors of $\tau$ are needed to label the ``mismatch'' in filtrations. However, we still need only a single factor, even though the difference in CE filtrations is now 2 instead of 1.

\printbibliography

\end{document}